\documentclass[aos,citesort,dvips]{arximspdf}
\usepackage{dcolumn}

\doi{10.1214/10-AOS831}
\volume{39}
\issue{2}
\pubyear{2011}
\firstpage{673}
\lastpage{701}

\makeatletter

\newcolumntype{d}[1]{D{.}{.}{#1}}

\newtheorem{theorem}{Theorem}
\newtheorem{lemma}{Lemma}

\newproclaim{definition}{Definition}
\newproclaim{example}{Example}

\newcommand{\IID}{i.i.d.}

\newcommand{\Natural}{\mathbb{N}}
\newcommand{\Rational}{\mathbb{Q}}
\newcommand{\Real}{\mathbb{R}}
\newcommand{\bsa}{\mathbf{a}}
\newcommand{\bsb}{\mathbf{b}}
\newcommand{\bsc}{\mathbf{c}}
\newcommand{\bsu}{\mathbf{u}}
\newcommand{\bsx}{\mathbf{x}}
\newcommand{\bsy}{\mathbf{y}}
\newcommand{\bsv}{\mathbf{v}}
\newcommand{\bszero}{\mathbf{0}}
\newcommand{\bsdelta}{\bolds{\delta}}
\newcommand{\bsone}{\mathbf{1}}
\newcommand{\bsZ}{\mathbf{Z}}
\newcommand{\Bcal}{\mathcal{B}}
\newcommand{\Ccal}{\mathcal{C}}
\newcommand{\Ncal}{\mathcal{N}}
\newcommand{\rd}{{d}}
\newcommand{\vol}{\operatorname{Vol}}
\newcommand{\wt}{\widetilde}

\newcommand{\hk}{\mathrm{HK}}
\newcommand{\xset}{\Omega}

\newcommand{\ds}{}
\newcommand{\diag}{\operatorname{diag}}
\newcommand{\tran}{\mathsf{T}}

\makeatother

\begin{document}
\begin{frontmatter}

\title{Consistency of Markov chain quasi-Monte Carlo\\ on continuous
state spaces}
\runtitle{QMC for MCMC}

\begin{aug}
\author[A]{\fnms{S.} \snm{Chen}\thanksref{t1}},
\author[B]{\fnms{J.} \snm{Dick}\thanksref{t2}} and
\author[C]{\fnms{A. B.} \snm{Owen}\corref{}\thanksref{t1}\ead[label=e1]{owen@stat.stanford.edu}}
\runauthor{S. Chen, J. Dick and A. B. Owen}
\affiliation{Stanford University, University of New South Wales and\\
Stanford University}
\address[A]{S. Chen\\
Stanford University\\
Sequoia Hall\\
Stanford, California 94305\\
USA} %adresu isvedimo komanda gale!
\address[B]{J. Dick\\
School of Mathematics and Statistics\\
University of New South Wales\\
Sydney\\
Australia}
\address[C]{A. B. Owen\\
Stanford University\\
Sequoia Hall\\
Stanford, California 94305\\
USA\\
\printead{e1}}
\end{aug}

\thankstext{t1}{Supported by NSF Grants
DMS-06-04939 and DMS-09-06056.}
\thankstext{t2}{Supported by a Queen Elizabeth II fellowship.}

% HISTORY:
\received{\smonth{8} \syear{2009}}
\revised{\smonth{5} \syear{2010}}

% ABSTRACT
%
\begin{abstract}
The random numbers driving Markov chain Monte Carlo (MCMC) simulation
are usually modeled as independent $U(0,1)$ random variables.
Tribble [Markov chain Monte Carlo algorithms using completely uniformly
distributed driving sequences (2007) Stanford Univ.] reports
substantial improvements when those
random numbers are replaced by carefully balanced inputs from
completely uniformly distributed sequences. The previous theoretical
justification for using anything other than \IID\  $U(0,1)$ points shows
consistency for estimated means, but only applies for discrete
stationary distributions. We extend those results to some MCMC
algorithms for continuous stationary distributions. The main motivation
is the search for quasi-Monte Carlo versions of MCMC. As a side
benefit, the results also establish consistency for the usual method of
using pseudo-random numbers in place of random ones.
\end{abstract}

% KEYWORDS
%

\end{frontmatter}

%s1 ###
\section{Introduction}

In Markov chain Monte Carlo (MCMC), one simulates
a Markov chain and uses sample averages to
estimate corresponding means of the stationary
distribution of the chain.
MCMC has become a staple tool in the physical sciences
and in Bayesian statistics.
When sampling the Markov chain, the transitions are driven by
a stream of independent $U(0,1)$ random numbers.

In this paper, we study what happens when the \IID\  $U(0,1)$
random numbers are replaced by deterministic sequences,
or by some dependent $U(0,1)$ values.
The motivation for replacing \IID\  $U(0,1)$ points is
that carefully stratified inputs
may lead to more accurate sample averages.
One must be cautious though, because as with adaptive MCMC
\cite{haarsakstamm2001,andrmoul2006},
the resulting simulated points do not have the Markov property.

The utmost in stratification is provided by
quasi-Monte Carlo (QMC) points.
There were a couple of attempts at merging
QMC into MCMC around 1970, and then again
starting in the late 1990s. It is only
recently that significant improvements have
been reported in numerical investigations.
For example, Tribble \cite{trib2007} reports variance
reductions of several thousand fold and an apparent
improved convergence rate for some Gibbs sampling problems.
Those results motivate our theoretical work. They are
described more fully in the literature survey below.

To describe our contribution,
represent MCMC sampling
via $\bsx_{i+1} = \phi(\bsx_i,\break\bsu_i)$ for $i=1,\ldots,n$,
where $\bsx_0$ is a nonrandom starting point
and $\bsu_i \in(0,1)^d$.
The points $\bsx_i$ belong to a state space $\Omega\subset\Real^s$.
The function $\phi$ is chosen so that $\bsx_i$ form
an ergodic Markov chain
with the desired stationary distribution $\pi$
when $\bsu_i\sim U(0,1)^d$ independently.
For a bounded continuous function $f\dvtx\Omega\to\Real$,
let $\theta(f) = \int_\Omega f(\bsx)\pi(\bsx)\,\rd\bsx$
and $\hat\theta_n(f) = (1/n)\sum_{i=1}^n f(\bsx_i)$.
%If $f$ is a bounded and continuous,
Then $\hat\theta_n(f)\to_{\mathbb{P}} \theta(f)$
%in probability
as $n\to\infty$.
In this paper, we supply sufficient conditions on $\phi$ and on the
deterministic sequences $\bsu_i$
so that $\hat\theta_n(f)\to\theta(f)$ holds
when those deterministic sequences are used instead of random ones.
The main condition is that the components of $\bsu_i$ be
taken from a completely uniformly distributed (CUD) sequence,
as described below.

Ours are the first results to prove
that deterministic sampling applied to MCMC problems on continuous state
spaces is consistent.
In practice, of course, floating point computations
take place on a large discrete state space. But invoking
finite precision does not provide a satisfying description of
continuous MCMC problems. In a finite state space
argument, the resulting state spaces are so big
that vanishingly few states will ever be visited
in a given simulation. Then if one switches from $32$ to $64$
to $128$ bit representations, the problem seemingly requires
vastly larger sample sizes, but in reality is not
materially more difficult.

To avoid using the finite state shortcut, we adopt a
computational model with infinite precision.
As a side benefit, this paper shows that the standard
practice of
replacing genuine \IID\  values $\bsu_i$ by deterministic
pseudo-random numbers is consistent for some problems
with continuous state spaces. We do not think
many people doubted this, but neither has it been established
before, to our knowledge.
It is already known from Roberts, Rosenthal and Schwartz \cite
{roberoseschw1998}
that, under certain conditions, a geometrically ergodic Markov
chain remains so under small perturbations, such as rounding.
That work does not address the replacement of random points by
deterministic ones
that we make here.

% Chentsov \cite{chen1967} proved consistency
% for some quasi-Monte Carlo approaches to sampling
% Markov chains, when the inverse CDF is used to
% generate each transition and inversion is also
% used to sample from an initial distribution.
% Owen and Tribble \cite{qmcmetro} extended the argument
% from inversion to Metropolis-Hastings and Gibbs sampling
% subject to a Riemann integrability condition on
% the transition mechanism.

% The proven mathematical properties of QMC points are related
% to equidistribution over rectangular regions. These are the
% same properties used to verify that pseudo-random number generators
% properly simulate randomness. A consequence of our work is that
% replacing genuine \IID\  random variables by pseudo-random numbers
% also leads to consistent MCMC estimates. We do not think anybody
%doubted
% this, but neither has it been previously proved for the
%Metropolis-Hastings
% algorithm outside of the finite state space setting.
%s1.1 ###
\subsection{Literature review}

There have been a small number of prior attempts to apply
QMC sampling to MCMC problems.
The first appears to have been
Chentsov \cite{chen1967}, whose work
appeared in 1967, followed
by Sobol'~\cite{sobo1974} in 1974.
Both papers assume that the Markov chain has
a discrete state space and that the transitions
are sampled by inversion.
Unfortunately, QMC does not usually bring large
performance improvements on such unsmooth problems
and inversion is not a very convenient method.

Chentsov replaces \IID\  samples by one long
CUD sequence, and this is the method we
will explain and then adapt to continuous problems.
Sobol' uses what is conceptually an $n\times\infty$
matrix of values from the unit interval. Each
row is used to make transitions until the chain returns to its
starting state. Then the sampling starts using
the next row.
It is like deterministic regenerative sampling.
Sobol' shows that the error converges as $O(1/n)$
in the very special case where the transition probabilities
are all rational numbers with denominator a power of $2$.
These methods were not widely cited and,
until recently, were almost forgotten, probably due
to the difficulty of gaining large improvements
in discrete problems, and the computational awkwardness of inversion
as a transition mechanism for discrete state spaces.

The next attempt that we found is that of Liao \cite{liao1998}
in 1998. Liao takes a~set of QMC points in $[0,1]^d$
shuffles them in random order, and uses them to drive
an MCMC. He reports 4- to 25-fold efficiency improvements, but gives
no theory. An analysis of Liao's method
is given in \cite{qmcmetro2}.
Later, Chau\-dary~\cite{chau2004} tried a different strategy
using QMC to generate balanced proposals for Metropolis--Hastings
sampling, but found only small improvements and did not
publish the work.
Craiu and Lemieux \cite{crailemi2007} also consider multiple-try
Metropolis and find variance reductions of up to $30$\%, which
is still modest. Earlier, Lemieux and Sidorsky \cite{lemisido2006}
report variance reduction factors ranging from about $1.5$
to about $18$ in some work using QMC in conjunction with
the perfect sampling method of Propp and Wilson \cite{propwils1996}.

Only recently
have there been significantly large benefits
from the combination of QMC and MCMC. Those benefits have
mainly arisen for problems on continuous state spaces.
Tribble's \cite{trib2007} best results come from
Gibbs sampling problems computing posterior means.
For problems with $d$ parameters, he used every
$d$-tuple from a small custom built
linear feedback shift register (LFSR).
One example is the well-known model used by Gelfand and Smith \cite
{gelfsmit1990}
for failure events of $10$ pumps from the article by Gaver and
O'Murcheartaigh \cite{gaveomui1987}.
There are $11$ unknown parameters, one for each pump and one for the scale
parameter in the distribution of pump failure rates.
A second example is a $42$ parameter probit model for vasorestriction
based on a famous data set from \cite{finn1947} and analyzed
using latent variables as in Albert and Chib \cite{albechib1993}.
Of those $42$ parameters, the $3$ regression coefficients are of
greatest interest and $39$ latent variables are nuisance variables.
Table \ref{tab:tribbles} sets out variance reduction factors
found for randomized CUD versus \IID\  sampling. The improvements
appear to grow with $n$, and are evident at very small sample sizes.

%t1 ###
%
\begin{table}
\caption{Variance reduction factors from Tribble \protect\cite{trib2007} for
two Gibbs sampling problems. For the pumps data, the greatest and least
variance reduction for a randomized CUD sequence versus \IID\  sampling is
shown. For the vasorestriction data, greatest and least variance
reductions for the three regression parameters are shown. See
\protect\cite{trib2007} for simulation details}\label{tab:tribbles}
\begin{tabular*}{\tablewidth}{@{\extracolsep{\fill}}ld{3.0}d{4.0}d{3.0}d{4.0}d{4.0}d{5.0}@{\hspace*{-2pt}}}
\hline
& \multicolumn{2}{c}{$\bolds{n=2^{10}}$} & \multicolumn{2}{c}{$\bolds{n=2^{12}}$}
& \multicolumn{2}{c@{}}{$\bolds{n=2^{14}}$}\\[-4pt]
& \multicolumn{2}{c}{\hrulefill} & \multicolumn{2}{c}{\hrulefill}
& \multicolumn{2}{c@{}}{\hrulefill}\\
\textbf{Data} & \multicolumn{1}{c}{\textbf{min}} & \multicolumn{1}{c}{\textbf{max}}
& \multicolumn{1}{c}{\textbf{min}} & \multicolumn{1}{c}{\textbf{max}} & \multicolumn{1}{c}{\textbf{min}}
& \multicolumn{1}{c@{}}{\textbf{max}}\\
\hline
Pumps & 286 & 1543 & 304 & 5003 & 1186 & 16089\\
Vasorestriction & 14 & 15 & 56 & 76 & 108 & 124\\
\hline
\end{tabular*}
\end{table}

There is another line of research in which large
improvements have been obtained by combining
QMC with MCMC. This is the array-RQMC method described
in L'Ecuyer, Lecot and Tuffin \cite{leculecotuff2008}
and other articles.
That method simulates numerous chains in parallel
using quasi-Monte Carlo to update all the chains.
It requires a complicated method to match the
update variables for each step to the various evolving chains.
This method
has achieved variance reductions of many thousand
fold on some problems from queuing and finance.
Very few properties have been established for it,
beyond the case of heat particles in one
dimension that was considered by Morokoff and
Caflisch~\cite{morocafl1993}.

Finally, Jim Propp's rotor-router method
is a form of deterministic Markov chain sampling.
It has brought large efficiency improvements
for some problems on a discrete state space
and has been shown to converge at better than
the Monte Carlo rate on some problems.
See, for example, Doerr and Friedrich \cite{doerfrie2009}.

The use of CUD sequences that we study has one
practical advantage
compared to the rotor-router, array-RQMC,
regenerative sampling, and the other methods.
It only requires replacing
the \IID\  sequence used in a typical MCMC run by
some other list of numbers.

%s1.2 ###
\subsection{Outline}

The paper is organized around our main results
which appear in Section \ref{sec:consistency}.
Theorem \ref{thm:withhomecouple}
gives sufficient conditions for consistency of QMC-MCMC
sampling by Metropolis--Hastings.
Theorem \ref{thm:ifmAlsFuh} gives sufficient conditions
for consistency of QMC-MCMC sampling for the systematic
scan Gibbs sampler.

Section \ref{sec:background} contains necessary background
and notation for the two main theorems of Section \ref{sec:consistency}.
It introduces quasi-Monte Carlo and Markov chain Monte Carlo
giving key definitions we need in each case.
That section presents the Rosenblatt--Chentsov transformation.
We have combined a classic sequential inversion method based
on the Rosenblatt transformation with an elegant coupling
argument that Chentsov \cite{chen1967} used.
%There is also one technical lemma that fits there.

The consistency results for Metropolis--Hastings
(Theorem \ref{thm:withhomecouple})
make moderately strong
assumptions in order to ensure that a coupling occurs.
Section~\ref{sec:homeandcouple}
shows that those assumptions are satisfied by some Metropolized
independence samplers and also by some slice samplers.
We also assumed some Riemann integrability properties
for our MCMC proposals. The Riemann integral is awkward
compared to the Lebesgue integral, but considering it is
necessary when we want to study specific algorithms on deterministic inputs.
Section \ref{sec:riemaninteg}
gives sufficient conditions for an MCMC algorithm to satisfy
the required Riemann integrability conditions.

Our consistency results for the Gibbs sampler
(Theorem \ref{thm:ifmAlsFuh}) require some contraction properties
and some Jordan measurability.
Section \ref{sec:gibbsexample} shows that these properties
hold under reasonable conditions.
Section \ref{sec:openclose} has a brief discussion on
open versus closed intervals for uniform random numbers.
Our conclusions are in Section \ref{sec:conclusions}.
The lengthier or more technical proofs are placed
in the \hyperref[app]{Appendix}.\vadjust{\goodbreak}

%s2 ###
\section{Background on QMC and MCMC}\label{sec:background}
%s2.1 ###
\subsection{Notation}

Our random vectors are denoted
by $\bsx= (x_1,\ldots,x_s)\in\xset\subseteq\Real^s$
for $s\ge1$.
%When $s=1$ the vectors become scalars and may
%be written either as $\bsx$ or $x$ depending on context.
Points in the unit cube $[0,1]^d$ are
denoted by $\bsu= (u_1,\ldots, u_d)$.
Two points $\bsa,\bsb\in\Real^d$ with $a_j<b_j$ for $j=1,\ldots,d$
define a rectangle $\prod_{j=1}^d [a_j,\break b_j]$, denoted
by $[\bsa,\bsb]$ for short.
The indicator (or characteristic)
function of a set $A\subset\Real^d$ is written $1_A$.
%$$1_A(\bsx) = \begin{cases} 1 & \bsx\in A\\0 &\text{else.}\end{cases}$$
%That is $1_A(\bsx) =1$ for $\bsx\in A$ and $1_A(\bsx)=0$ for $\bsx\not

We assume the reader is familiar with the definition
of the (proper) Riemann integral, for a bounded function
on a finite rectangle $[\bsa,\bsb]\subset\Real^d$.
%We use $\lambda_d$ to denote Lebesgue measure on $\Real^d$.
The bounded set $A\subset\Real^d$ is Jordan measurable if
%the function
%$$1_A(\bsx) = \begin{cases} 1 & \bsx\in A\\0 &\text{else}
%its indicator (or characteristic) function
$1_A$ is Riemann integrable on a~bounded rectangle
containing $A$. By Lebesgue's theorem
(see Section \ref{sec:riemaninteg})~$A$ is Jordan measurable
if $\lambda_d(\partial A) =0$. Here $\lambda_d$
denotes Lebesgue measure on $\Real^d$, and $\partial A$ is
the boundary of $A$, that is, the set on which
$1_A$ is discontinuous.

%s2.2 ###
\subsection{QMC background}

Here, we give a short summary of quasi-Monte Carlo.
Further information may be found in the monograph
by Niedereiter~\cite{nied92}.

QMC is ordinarily used to approximate integrals
over the unit cube $[0,1]^d$, for $d\in\Natural$.
Let $\bsx_1,\ldots,\bsx_n\in[0,1]^d$.
The QMC estimate of $\theta(f) = \int_{[0,1]^d} f(\bsx)\,\rd\bsx$ is
$\hat\theta_n(f) = \frac1n\sum_{i=1}^nf(\bsx_i)$,
just as we would use in plain Monte Carlo.
The difference is that in QMC, distinct points $\bsx_i$ are chosen
%%% this is very picky. If the points are not distinct
%%% then the probability on x_i is not 1/n. It is the
%%% n_i/n where n_i is the number of points with x_j=x_i
deterministically to make the discrete probability distribution
with an atom of size $1/n$ at each $\bsx_i$
close to the continuous $U[0,1]^d$ distribution.

The distance between these distributions is quantified
by discrepancy measures.
The local discrepancy of $\bsx_1,\ldots,\bsx_n$ at $\bsa\in[0,1]^d$ is
%
%e1 ###
%
\begin{equation}\label{eq:localdis}
\delta(\bsa)=\delta(\bsa;\bsx_1,\ldots,\bsx_n)=
\frac1n\sum_{i=1}^n1_{[\bszero,\bsa)}(\bsx_i)
-\prod_{j=1}^da_j.
\end{equation}
The star discrepancy of $\bsx_1,\ldots,\bsx_n$ in dimension $d$ is
%
%e2 ###
%
\begin{equation}\label{eq:stardisc}
D_n^{*d} =D_n^{*d}(\bsx_1,\ldots,\bsx_n) ={\sup_{\bsa\in
[0,1]^d}}|\delta(\bsa;\bsx_1,\ldots,\bsx_n)|.
\end{equation}
For $d=1$, the star discrepancy
reduces to the Kolmogorov--Smirnov distance
between a discrete and a continuous uniform distribution.

A uniformly distributed sequence is one for which $D_n^{*d}\to0$
as $n\to\infty$.
If $\bsx_i$ are uniformly distributed then $\hat\theta_n(f)\to
\theta(f)$
provided that $f$ is Riemann integrable.

Under stronger conditions than Riemann integrability,
we can get rates of convergence for QMC.
The Koksma--Hlawka inequality is
%
%e3 ###
%
\begin{equation}\label{eq:khbound}
|\hat\theta_n(f)-\theta(f)| \le D_n^{*d} V_{\hk}(f),
\end{equation}
where $V_{\hk}$ is the total variation of $f$ in
the sense of Hardy and Krause. For properties of $V_{\hk}$
and other multidimensional variation measures,
see \cite{variation}.

Equation (\ref{eq:khbound}) gives a deterministic upper bound
on the integration error, and it factors into a measure
of the points' quality and a measure of the integrand's roughness.
There exist constructions $\bsx_1,\ldots,\bsx_n$ where
$D_n^{*d} = O(n^{-1+\epsilon})$ holds for any $\epsilon>0$.
Therefore, functions of finite variation can be integrated
at a much better rate by QMC than by MC. Rates of convergence of
$O(n^{-\alpha} (\log n)^{d \alpha})$, where $\alpha\ge1$ denotes
the smoothness of the integrand which can therefore be arbitrarily
large, can also be achieved~\cite{dickmcqmc2009}.

%%% Josef: Here we could put in a sentence pointing to a survey
% of higher accuracy QMC methods, such as might be in your
% mcqmc 2008 article with Baldeaux.

%There are many other discrepancies, using more general sets than
%$[\bszero,\bsa]$ and also replacing the supremum by an $L^p$ norm.

Equation (\ref{eq:khbound}) is not usable for error
estimation. Computing the star discrepancy is very difficult
\cite{gnewsrivwinz2008},
and computing $V_\hk(f)$ is harder than integrating $f$.
Practical error estimates for QMC may be obtained using randomized quasi-Monte
Carlo (RQMC). In RQMC each $\bsx_i\sim U[0,1]^d$ individually
while the ensemble $\bsx_1,\ldots,\bsx_n$ has
$\Pr( D_n^{*d}(\bsx_1,\ldots,\bsx_n)<C(\log n)^d/n)=1$
for some $C<\infty$.
For an example, see \cite{rtms}.
A small number of independent replicates of the RQMC estimate
can be used to get an error estimate. RQMC has the further
benefit of making QMC unbiased.
For a survey of RQMC, see \cite{leculemi1999a}.

A key distinction between QMC and MC is that the
former is effective for Riemann integrable functions,
while the latter, in principle, works for Lebesgue
integrable functions. In practice, MC is usually implemented
with deterministic pseudo-random numbers. The best generators
are proved to simulate independent $U[0,1]$ random variables
based on either discrepancy measures over rectangles or on spectral
measures. Those conditions are enough to prove convergence
for averages of Riemann integrable functions, but not for
Lebesgue integrable functions.
As a result, ordinary Monte Carlo with pseudo-random numbers
is also problematic for Lebesgue integrable functions that
are not Riemann integrable.

%s2.3 ###
\subsection{Completely uniformly distributed}

In the Markov chain context, we need a lesser known QMC concept as follows.
%The consecutive $d$-tuples from a random number generator
%should be nearly uniformly distributed. Ideally this holds
%for all $d\ge1$.
A sequence $u_1,u_2,\ldots\in[0,1]$ is completely uniformly
distributed (CUD)
if for any $d\ge1$ the points
$\bsx_i^{(d)}=(u_i,\ldots,u_{i+d-1})$ satisfy
$D_n^{*d}(\bsx^{(d)}_1,\ldots,\bsx^{(d)}_n)\to0$ as $n\to\infty$.
This is one of the definitions of a random sequence from
Knuth \cite{knut199723}, and it is an important property
for modern random number generators.

Using a CUD sequence in an MCMC is akin to using
up the entire period of a random number generator,
as remarked by Niederreiter \cite{nied1986} in 1986. It is
then necessary to use a small random number generator.
The CUD sequences used by Tribble \cite{trib2007} are
miniature versions of linear congruential generators
and feedback shift register generators. As such, they
are no slower than ordinary pseudo-random numbers.

In the QMC context, we need to consider nonoverlapping
$d$-tuples $\wt\bsx_i^{(d)}=(u_{di-d+1},\ldots,u_{di})$ for $i\ge1$.
It is known \cite{chen1967} that
%
%e4 ###
%
\begin{eqnarray}\label{eq:overlapornot}
&D_n^{*d}\bigl(\bsx^{(d)}_1,\ldots,\bsx^{(d)}_n\bigr)\to0\qquad \forall d\ge1,&
\nonumber\\
&\iff&\\
&D_n^{*d}\bigl(\wt\bsx^{(d)}_1,\ldots,\wt\bsx^{(d)}_n\bigr)\to0\qquad
\forall d\ge
1.&\nonumber
\end{eqnarray}

%s2.4 ###
\subsection{MCMC iterations}

In the QMC context, the function $f$ subsumes
all the necessary transformations to turn a finite
list of \IID\  $U[0,1]$ random variables into
the desired nonuniformly distributed quantities, as well as the
function of those quantities whose expectation
we seek.
In some problems, we are unable to find such transformations,
and so we turn to MCMC methods.

Suppose that we want to sample $\bsx\sim\pi$
for a density function $\pi$ defined with
respect to Lebesgue measure on $\xset\subseteq\Real^s$.
For definiteness, we will seek to approximate
$\theta(f) = \int_\xset f(\bsx)\pi(\bsx)\,\rd\bsx$.
In this section, we briefly present MCMC.
For a full description of MCMC, see the monographs by
Liu \cite{liu2001} or Robert and Casella \cite{robecase2004}.

In an MCMC simulation,
we choose an arbitrary $\bsx_0\in\xset$ with $\pi(\bsx_0)>0$
and then for $i\ge1$ update via
%
%e5 ###
%
\begin{equation}\label{eq:defphi}
\bsx_{i} = \phi(\bsx_{i-1},\bsu_{i}),
\end{equation}
where $\bsu_i\in[0,1]^d$ and $\phi$
is an update function described below.
The distribution of $\bsx_i$
depends on $\bsx_0,\ldots,\bsx_{i-1}$ only through $\bsx_{i-1}$
and so these random variables have the Markov property.
The function $\phi$ is chosen so that the stationary
distribution of $\bsx_i$ is $\pi$.
Then we estimate $\theta(f)$ by
$\hat\theta_n(f)=\frac1n\sum_{i=1}^nf(\bsx_i)$ as before.
If a burn-in period was used, we assume that $\bsx_0$
is the last point of it.

First, we describe the Metropolis--Hastings algorithm
for computing $\phi(\bsx,\bsu)$ from the current
point $\bsx\in\xset$ and $\bsu\in[0,1]^d$.
It begins with a proposal $\bsy$ taken
from a transition kernel $P(\bsx, \rd\bsy)$.
With genuinely random proposals, the transition kernel
gives a complete description. But for either quasi-Monte
Carlo or pseudo-random sampling, it matters how we
actually generate the proposal. We will assume that $d-1$
$U[0,1]$ random variables are used to generate $\bsy$
via $\bsy= \psi_{\bsx}(u_{1\dvtx(d-1)})$.
Then the proposal $\bsy$ is either accepted or rejected with
probability $A(\bsx,\bsy)$. The decision is typically based
on whether the $d$th random variable $u_d$ is
below $A$.
\begin{definition}[(Generator)]
The function $\psi\dvtx[0,1]^d\to\Real^s$ is a generator
for the distribution $F$ on $\Real^s$ if
$\psi(\bsu)\sim F$ when $\bsu\sim U[0,1]^d$.
\end{definition}
\begin{definition}[(Metropolis--Hastings update)]
For $\bsx\in\xset$, let $\psi_{\bsx}\dvtx[0,1]^{d-1}\to\xset$
be a generator for the transition kernel $P(\bsx, \rd\bsy)$
with conditional density \mbox{$p(\cdot\mid\bsx)$}.
The Metropolis--Hastings sampler has
\[
\phi(\bsx,\bsu)
=
\cases{
\bsy(\bsx,\bsu), &\quad $u_d \le A(\bsx,\bsu)$,\cr
\bsx, &\quad $u_d > A(\bsx,\bsu)$,}
\]
where
$\bsy(\bsx,\bsu) = \psi_{\bsx}(\bsu_{1\dvtx(d-1)})$
and
\[
A(\bsx,\bsu) = \min\biggl(1,\frac{\pi(\bsy(\bsx,\bsu)) p(\bsx
\mid{\bsy(\bsx,\bsu))}
}
{\pi(\bsx) p(\bsy(\bsx,\bsu)\mid\bsx)}\biggr).
\]
\end{definition}
\begin{example}[{[Metropolized independence sampler (MIS)]}]
The MIS update is a special
case of the Metropolis--Hastings update in
which $\bsy(\bsx,\bsu) = \psi(\bsu_{1\dvtx(d-1)})$ does
not depend on $\bsx$.
%Instead $\psi$ is a generator for a proposal distribution $F$.
\end{example}
\begin{example}[{[Random walk Metropolis (RWM)]}]
The RWM update is a special case of the Metropolis--Hastings
update in which $\bsy(\bsx,\bsu)=\bsx+\psi(\bsu_{1\dvtx(d-1)})$
for some generator $\psi$ not depending on $\bsx$.
\end{example}
\begin{definition}[(Systematic scan Gibbs sampler)]
Let $\bsx= (x_1,\ldots, x_s) \in\Real^d$ with $x_j \in\Real
^{k_j}$ and $d = \sum_{j=1}^s k_j$.
To construct the systematic scan Gibbs sampler,
let $\psi_{j,\bsx_{-j}}(\bsu_j)$ be a $k_j$-dimensional
generator of the full conditional distribution
of $x_j$ given $x_\ell$ for all $\ell\ne j$.
This Gibbs sampler generates the new point
using $\bsu\in[0,1]^d$.
Write $\bsu= (\bsu_1,\ldots,\bsu_s)$ with $\bsu_j\in[0,1]^{k_j}$.
%% we may need notation for vector of vectors, tuples, row vs col etc
%% the less to say, the better, since it is obvious
The systematic scan Gibbs sampler has
\[
\phi(\bsx,\bsu) =
(\phi_1(\bsx,\bsu),\phi_2(\bsx,\bsu),\ldots,\phi_s(\bsx,\bsu)),
\]
where, for $1\le j\le s$,
\[
\phi_j(\bsx,\bsu) = \psi_{j, \bsx_{[j]}}(\bsu_j)
\]
and
$\bsx_{[j]} = (\phi_1(\bsx,\bsu),\ldots,\phi_{j-1}(\bsx,\bsu),
x_{j+1},\ldots,x_d)$.
%This is a $d$--dimensional generator where $d=\sum_{j=1}^sk_j$.
\end{definition}
\begin{example}[(Inversive slice sampler)]
Let $\pi$ be a probability density function on
$\xset\subset\Real^s$.
Let $\xset' =
\{ (y,\bsx)\mid\bsx\in\xset, 0\le y\le\pi(\bsx)\}
\subset\Real^{s+1}$
and let $\pi'$ be the uniform distribution on
$\xset'$.
The inversive slice sampler is the systematic scan Gibbs
sampler for $\pi'$ with each $k_j=1$
using inversion for every $\psi_{j,\bsx_{[j]}}$.
\end{example}

There are many other slice samplers.
See \cite{neal2003}.
It is elementary that $(y,\bsx)\sim\pi'$
implies $\bsx\sim\pi$.
It is more usual to use $(\bsx,y)$, but
our setting simplifies when we assume $y$
is updated first.

%s2.5 ###
\subsection{Some specific generators}

We generate our random variables as functions
of independent uniform random variables.
The generators we consider require a finite
number of inputs, so acceptance-rejection is
not directly covered, but see the note in Section \ref{sec:conclusions}.

For an encyclopedic presentation of methods
to generate nonuniform random vectors,
see Devroye \cite{devr1986}.
Here, we limit ourselves to inversion and some generalizations
culminating in the Rosenblatt--Chentsov transformation introduced below.
We will not need to assume that $\pi$ can be sampled by inversion.
We only need inversion for an oracle used later in a coupling argument.

Let $F$ be the CDF of $x\in\Real$, and
for $0< u <1$ define
\[
F^{-1}(u) = \inf\{ x\mid F(x)\ge u\}.
\]
Take $F^{-1}(0)=\lim_{u\to0^+}F^{-1}(u)$
and $F^{-1}(1)=\lim_{u\to1^-}F^{-1}(u)$, using extended
reals if necessary.
Then $x = F^{-1}(u)$ has distribution
$F$ on $\Real$ when $u\sim U[0,1]$.

Multidimensional inversion is based on inverting
the Rosenblatt transformation~\cite{rose1952}.
Let $F$ be the joint distribution of $\bsx\in\Real^s$.
Let $F_1$ be the marginal CDF of $x_1$
and for $j=2,\ldots,s$, let $F_j( \cdot;\bsx_{1\dvtx(j-1)})$ be the
conditional CDF of $x_j$ given $x_1,\ldots,x_{j-1}$.
The inverse Rosenblatt transformation $\psi_R$
of $\bsu\in[0,1]^s$ is $\psi_R(\bsu)=\bsx\in\Real^s$
where
\[
x_1 = F_1^{-1}(u_1)
\]
and
\[
x_j = F_j^{-1}\bigl(u_j;\bsx_{1\dvtx(j-1)}\bigr),\qquad j\ge1.
\]
If $\bsu\sim U[0,1]^s$, then $\psi_R(\bsu)\sim F$.

%In the Rosenblatt transformation,
%the dimension $s$ of $\bsx$ matches that of $\bsu$.
%In many other transformations,
%a $d$--dimensional random vector is constructed
%from $s$ uniform random variables where $s>d$.

We will use the inverse Rosenblatt transformation as a first
step in a~coupling argument which extends
the one in Chentsov \cite{chen1967}.
\begin{definition}[(Rosenblatt--Chentsov transformation)]
Let $\psi_R$ be the inverse Rosenblatt transformation
for the stationary distribution $\pi$
and let $\phi$ be the update function for MCMC.
The Rosenblatt--Chentsov transformation of
the finite sequence $\bsu_0,\bsu_1,\ldots,\bsu_m\in[0,1]^d$
is the finite sequence $\bsx_0,\ldots, \bsx_m\in\xset\subset\Real^s$,
with $s\le d$,
where $\bsx_0=\psi_R(\bsu_{0,1\dvtx s})$
and $\bsx_i=\phi(\bsx_0,\bsu_i)$ for $i=1,\ldots,m$.
\end{definition}

%% xxx we will need to impose $d\ge s$ I think

The Rosenblatt--Chentsov transformation starts off using
$\bsu_0$ and inversion to generate $\bsx_0$ and then
it applies whatever generators are embedded in $\phi$
with the innovations $\bsu_i$,
to sample the transition kernel. The transition function $\phi$ need not
be based on inversion.

%s3 ###
\section{Consistency for MCQMC sampling}\label{sec:consistency}

In this section, we prove sufficient conditions for some
deterministic MCQMC samplers to sample consistently.
The same proof applies to deterministic pseudo-random sampling.
First, we define consistency, then some
regularity conditions, and then we give the main results.

%s3.1 ###
\subsection{Definition of consistency}

Our definition of consistency is that the empirical distribution
of the MCMC samples converges weakly to $\pi$.
\begin{definition}\label{def:consist}
The triangular array $\bsx_{n,1},\ldots,\bsx_{n,n}\in\Real^s$
for $n$ in an infinite set $\Natural^*\subset\Natural$ consistently
samples the probability density function $\pi$~if
%$\hat\pi_n$ converges weakly to $\pi$ where
%$\hat\pi_n$ is $1/n$ times the sum of point
%mass distributions at $\bsx_{n,i}$ for $i=1,\ldots,n$.
%= \int_{[\bsa,\bsb]} \pi(\bsx)\rd\bsx
%holds for all rectangles $[\bsa,\bsb]\subseteq\Real^s$
%of finite volume.
%
%e6 ###
%
\begin{equation}\label{eq:smoothto}
\mathop{\lim_{n\to\infty}}_{n\in\Natural^*}
\frac1n\sum_{i=1}^n f(\bsx_{n,i})
= \int f(\bsx) \pi(\bsx)\,\rd\bsx
\end{equation}
holds for all bounded continuous functions $f\dvtx\xset\to\Real$.
The infinite sequence $\bsx_1,\bsx_2,\ldots\in\Real^s$
consistently samples $\pi$ if
the triangular array of initial subsequences
with $\bsx_{n,i}=\bsx_i$ for $i=1,\ldots,n$ does.
\end{definition}

In practice, we use a finite list of vectors
and so the triangular array formulation is a
closer description of what we do.
However, to simplify the presentation and
avoid giving two versions of everything, we will
work only with the infinite sequence version
of consistency.
Triangular array versions of CUD sampling for discrete state spaces
are given in \cite{qmcmetro2}.

It suffices to use functions $f$ in a convergence-determining class.
For example, we may suppose that $f$ is
uniformly continuous \cite{ash1972}, % page 198
or that $f= 1_{(\bsa,\bsb]}$~\cite{bill1999}. % p18
When $\pi$ is a continuous distribution, we
may use $f= 1_{[\bsa,\bsb]}$.

%s3.2 ###
\subsection{Regularity conditions}

Here, we define some assumptions that we~need
to make on the MCMC update functions.
%
% \begin{definition}\label{def:homeonly}
% Let $\xset\subseteq\Real^s$ be a continuous state space
% with update $\phi:\xset\times[0,1]^d\to\xset$.
% The set $\Hcal\subseteq\xset$ is called a home region of $\phi$
% if there exists a Jordan measurable
% set $\Bcal\subseteq[0,1]^d$ of positive
% volume such that $\phi(\bsx,\bsu)\in\Hcal$
% whenever $\bsu\in\Bcal$.
% \end{definition}
%
% Wherever the walk is at the moment, it can return
% to the home region. Moreover when $\bsu\in\Bcal$
% then the walk goes home at the next step, no matter
% where it was.
%
% \begin{definition}\label{def:couplingreg}
% The region $\Hcal\subseteq\xset$ is called a
% coupling region if $\bsx,\bsx'\in\Hcal$
% implies that $\phi(\bsx,\bsu)=\phi(\bsx',\bsu)$
% for all $\bsu\in[0,1]^d$.
% \end{definition}
%
\begin{definition}\label{def:coupleregion}
Let $\Ccal\subset[0,1]^d$ have positive Jordan measure.
If $u\in\Ccal$ implies that
$\phi(\bsx,\bsu)=\phi(\bsx',\bsu)$ for all $\bsx,\bsx'\in
\Omega$,
then $\Ccal$ is a coupling region.
\end{definition}

Consider two iterations $\bsx_i=\phi(\bsx_{i-1},\bsu_i)$
and $\bsx'_i=\phi(\bsx'_{i-1},\bsu_i)$ with the
same innovations $\bsu_i$ but possibly different
starting points $\bsx_0$ and $\bsx'_0$.
If $\bsu_i\in\Ccal$, then
$\bsx_{j}=\bsx'_{j}$ holds for all $j \ge i$.
In Section \ref{sec:homeandcouple},
we give some nontrivial examples
of MCMC updates with coupling regions.
\begin{definition}[(Regular MCMC)]\label{def:regular}
Let $\bsx_m=\bsx_m(\bsu_0,\ldots,\bsu_m)$ be
the last point generated in the
Rosenblatt--Chentsov transformation,
viewed as a~function on $[0,1]^{d(m+1)}$.
%The MCMC is \textit{regular for rectangles}, if
%the set
%$\{ (\bsu_0,\ldots,\bsu_m)\mid\bsx_m(\bsu_0,\ldots,\bsu_m)\in[\bsa,\bsb]
%measurable whenever $[\bsa,\bsb]\subset\Real^s$ is bounded.
The MCMC is \textit{regular} (\textit{for bounded continuous functions}) if
the function $f(\bsx_m(\bsu_0,\ldots,\bsu_m))$ is Riemann
integrable on $[0,1]^{d(m+1)}$ whenever $f$ is bounded and continuous.
%The MCMC is regular if it satisfies either of these conditions.
%and is on a state space $\Omega\in\Real^s$ with $\lambda_s(\partial
% it may be that regularity implies Omega is well behaved
\end{definition}

Note that if an MCMC is regular, then the definition of the
Rosenblatt--Chentsov transformation implies that
\[
\int_{[0,1]^{d(m+1)}} f(\bsx_m(\bsu_0,\ldots,\bsu_m)) \,\rd\bsu_0
\cdots\rd\bsu_m =\int_{\xset} f(\bsx) \pi(\bsx) \,\rd\bsx
\]
for any $m \ge0$ and all bounded continuous functions $f$.

We can, of course, define regularity for MCMC also with respect to
other classes of functions. Indeed, there are numerous equivalent
conditions for regularity. For example, the Portmanteau
theorem (\cite{bill1999}, Chapter 1.2) implies that it is enough to
assume that the functions $f$ are bounded and uniformly continuous.
Of interest are also indicator functions of rectangles since they
appear in the definition of the local discrepancy at
(\ref{eq:localdis}). The following theorem states some
equivalent conditions. To simplify the statements, we write that
MCMC is \textit{regular for indicator functions} whenever
$1_{A}(\bsx_m(\bsu_0,\ldots,\bsu_m))$ is Riemann integrable on
$[0,1]^{d(m+1)}$, where $A$ is either $A = [\bsa,\bsb]$ with
$\bsa,\bsb$ finite or $A = \xset$.
\begin{theorem}\label{thm:eq_regularity}
The following statements are equivalent:
\begin{enumerate}[(iii)]
\item[(i)] MCMC is regular for bounded continuous functions.
\item[(ii)] MCMC is regular for bounded uniformly continuous functions.
\item[(iii)] MCMC is regular for indicator functions $1_{[\bsa,\bsb
]}$ of rectangles $[\bsa,\bsb]$.
\end{enumerate}
\end{theorem}
\begin{pf}
This result follows by applying the Portmanteau theorem
(\cite{bill1999}, Chapter~1.2)
and some methods from real analysis.
\end{pf}

A regular MCMC is one that satisfies any (and hence all) of the above.
%In the following we write for short regular by which we mean any of
%the equivalent statements: regular for bounded continuous functions,
%regular for bounded uniformly continuous functions, or regular for
%indicator functions of rectangles.

%s3.3 ###
\subsection{Main results for
Metropolis--Hastings}\label{subsec:MetropolisHastings}

Theorem \ref{thm:withhomecouple}
below is the main
result that we will use for Metropolis--Hastings sampling.
One does not expect CUD sampling to correct for an
MCMC algorithm that would not be ergodic
when sampled with \IID\  inputs.
Ergodicity is assured through our
assumption that there is %a home region that is also
a coupling region.
Section \ref{sec:homeandcouple} below shows that
some nontrivial MCMC methods have such regions.
Theorem \ref{thm:withhomecouple} does not require the detailed balance
condition that Metropolis--Hastings satisfies, and so it
may apply to some nonreversible chains too.
\begin{theorem}\label{thm:withhomecouple}
Let $\Omega\subseteq\Real^s$ and let $\bsx_0 \in\xset$, and for
$i \ge1$ let $\bsx_{i}=\phi(\bsx_{i-1},\bsu_i)$ where $\phi$ is the
update function of a regular MCMC with a coupling region~$\Ccal$.
If $\bsu_i = (v_{d(i-1)+1},\ldots,v_{di})$
for a CUD sequence $(v_i)_{i \ge1}$,
then $\bsx_1,\ldots,\bsx_n$ consistently samples $\pi$.
\end{theorem}

The proof of Theorem \ref{thm:withhomecouple} is in
the \hyperref[app]{Appendix}. It shows that the fraction of points $\bsx_i$
in a bounded rectangle $[\bsa,\bsb]$
converges to $\int_{[\bsa,\bsb]}\pi(\bsx)\,\rd\bsx$.
Almost the same proof technique applies to
expectations of bounded continuous functions.

%s3.4 ###
\subsection{Main results for Gibbs sampling}\label{subsec:Gibbs}

The Gibbs sampler can be viewed as a special
case of Metropolis--Hastings with acceptance
probability one.
However, it is more straightforward to study
it by applying results on iterated function
mappings to (\ref{eq:defphi})
using methods from Diaconis and Freedman~\cite{diacfree1999}
and Alsmeyer and Fuh \cite{alsmfuh2001}.

In this subsection, we assume that $(\xset,d)$ is a complete
separable metric space. We assume that the update function
$\phi(\bsx,\bsu)$ is jointly measurable in $\bsx$ and $\bsu$ and
that it is Lipschitz continuous in $\bsx$ for any $\bsu$. Lipschitz
continuity is defined through the metric $d(\cdot,\cdot)$ on $\xset
$. The
Lipschitz constant, which depends on $\bsu$, is
%
%e7 ###
%
\begin{equation}\label{eq:lips}
\ell(\bsu) = \sup_{\bsx\ne\bsx'}
\frac{d(\phi(\bsx,\bsu),\phi(\bsx',\bsu))}{d(\bsx
,\bsx')}.
\end{equation}
For each $\bsu_n\in[0,1]^d$, define
$L_n= \ell(\bsu_n)$.

Next, we present a theorem from Alsmeyer and Fuh \cite{alsmfuh2001}
on iterated random mappings.
The $n$ step iteration, denoted $\phi_n$,
is defined by
$\phi_1(\bsx;\bsu_1)=\phi(\bsx,\bsu_1)$
and for $n\ge2\dvtx
\phi_n(\bsx;\bsu_1,\ldots,\bsu_n)=\phi(\phi_{n-1}(\bsx;\bsu
_1,\ldots,\bsu_{n-1}),\bsu_n)$.
\begin{theorem}\label{thm:ifmAlsFuh}
Let the update function $\phi(\bsx,\bsu)$ be jointly measurable in
$\bsx$ and $\bsu$ with $\int_{[0,1]^d} \log(\ell(\bsu))\,\rd\bsu<0$
and, for some $p>0$, $\int_{[0,1]^d} \ell(\bsu)^p\,\rd\bsu<\infty$. Assume
that there is a point $\bsx'\in\xset$ with $\int_{[0,1]^d} \log^+(
d(\phi(\bsx',\bsu),\bsx'))\,\rd\bsu<\infty$ and $E(
d(\phi(\bsx',\bsu),\bsx')^p)<\infty$.
Then there is a $\gamma^\ast
\in(0,1)$ such that for all $\gamma\in(\gamma^\ast,1)$ there is a
$\alpha_\gamma\in(0,1)$ such that
for every $\bsx,\widehat{\bsx} \in\Omega$
%
%e8 ###
%
\begin{equation}\label{eq:ifm}
\lim_{m\to\infty} \alpha_\gamma^{-m}
\Pr\bigl( d(\phi_m(\bsx;\cdot),\phi_m(\widehat{\bsx};\cdot))
>\gamma^m\bigr)=0.
\end{equation}
%
%Then there is a $\gamma^\ast
%$\alpha_\gamma\in(0,1)$ such that
%>\gamma^m)=0,
%for all $\bsx,\widehat{\bsx} \in\Omega$.
\end{theorem}
\begin{pf}
This follows by specializing Corollary 2.5(a) of
\cite{alsmfuh2001} to the present setting.
\end{pf}
\begin{theorem}\label{thm:withiterfunction}
Let $(\Omega,d)$ be a complete separable metric space and let
$(v_i)_{i\ge1}$ be a CUD sequence such that for every sequence
$(d_n)_{n\ge1}$ of natural numbers with $d_n = O(\log n)$, we have
$\lim_{n\rightarrow\infty} D^{\ast d_n}_{n} = 0$. Let $\bsx_0 \in
\xset$, and for $i \ge1$ let $\bsx_{i}=\phi(\bsx_{i-1},\bsu_i)$ be
the Gibbs sampler update for stationary distribution $\pi$. Assume
that $\phi$ satisfies the conditions of Theorem \ref{thm:ifmAlsFuh}
and that there is a $\gamma\in(\gamma^\ast,1)$ such that
\[
\Bcal_{m}(\bsx,\widehat{\bsx}) = \{\bsv\in[0,1]^{dm}\dvtx
d(\phi_m(\bsx,\bsv), \phi_m(\widehat{\bsx},\bsv)) > \gamma^m \}
\]
is Jordan measurable for all $m \ge1$ and $\bsx,\widehat{\bsx} \in
\Omega$.
Under these conditions, if the Gibbs sampler is regular,
then $\bsx_1,\ldots, \bsx_n$ consistently samples $\pi$.
\end{theorem}

The proof of Theorem \ref{thm:withiterfunction}
is in the \hyperref[app]{Appendix}.
Like Theorem \ref{thm:withhomecouple}, it shows
that bounded rectangles $[\bsa,\bsb]$ have asymptotically
the correct proportion of points. Once again, similar
arguments apply for bounded continuous functions of $\bsx$.

Although not explicitly stated there, the proof of
\cite{dicknote2007}, Theorem 1, shows the existence of sequences
$(v_i)_{i\ge1}$ for which
\[
D^{\ast d}_{n}\bigl(\bigl\{\bigl(v_{d(i-1)+1}, \ldots, v_{di}\bigr), i=1,\ldots, n\bigr\}\bigr) \le
C\sqrt{\frac{d \log(n+1)}{n}},
\]
for all $n, d \in\Natural$, where $C > 0$ is a constant independent
of $n$ and $d$. Unfortunately, no explicit construction of such a
sequence is given in \cite{dicknote2007}. Then for any sequence
$(d_n)_{n \ge1}$ of natural numbers with $d_n = O(\log n)$ we
obtain that
\[
D^{\ast d_n}_{n}\bigl(\bigl\{\bigl(v_{d_n(i-1)+1}, \ldots, v_{d_n i}\bigr), i=1,\ldots,
n\bigr\}\bigr) \le C' \frac{ \log(n+1)}{\sqrt{n}} \rightarrow0
\qquad\mbox{as } n \rightarrow\infty.
\]

In Theorem \ref{thm:withhomecouple}, we assumed that the coupling
region $\Ccal$ is Jordan measurable
In Theorem \ref{thm:withiterfunction}, we do not have a coupling
region, but still have an analogous assumption, namely that the sets
$\Bcal_m(\bsx,\widehat{\bsx})$ are Jordan measurable. A condition on
$\phi$ which guarantees that $\Bcal_m(\bsx,\widehat{\bsx})$ is
Jordan measurable is given in Section \ref{sec:gibbsexample}.

%s4 ###
\section{Examples of coupling regions}\label{sec:homeandcouple}

Theorem \ref{thm:withhomecouple} used
coupling regions. These are somewhat special.
But they do exist for some realistic MCMC algorithms.
\begin{lemma}
Let $\phi$ be the update for the
Metropolized independence sampler on $\xset\subseteq\Real^s$
obtaining the proposal $\bsy=\psi(\bsu_{1\dvtx(d-1)})$,
where $\psi$ generates samples from the density $p$,
which are accepted when
\[
u_d \le\frac{\pi(\bsy)p(\bsx)}{\pi(\bsx)p(\bsy)}.
\]
Assume that the importance ratio is bounded above, that is,
\[
\kappa\equiv\sup_{\bsx\in\xset} \frac{\pi(\bsx)}{p(\bsx)}
<\infty.
\]
Suppose also that there is a rectangle $[\bsa,\bsb]\subset[0,1]^{d-1}$
of positive volume with
\[
\eta\equiv
\inf_{\bsu\in[\bsa,\bsb]} \frac{\pi(\psi(\bsu))}{p(\psi
(\bsu))}
>0.
\]
Then $\Ccal=[\bsa,\bsb]\times[0,\eta/\kappa]$ is a coupling region.
\end{lemma}
\begin{pf}
The set $\Ccal$ has positive Jordan measure.
%Let $\Hcal= \{\psi(\bsu_{1:(d-1)})\mid\bsu_{1:(d-1)}\in[\bsa,\bsb]\}$.
Suppose that $\bsu\in\Ccal$.~Then
\[
\pi(\bsy)p(\bsx)
\ge\eta p(\bsy) \frac1\kappa\pi(\bsx)\ge u_d p(\bsy)\pi(\bsx),
\]
and so $\phi(\bsx,\bsu)=\bsy$, regardless of $\bsx$.
\end{pf}
%
%CUD sampling of the independence sampler
%is consistent when
% $\sup_{\bsx\in\xset}\pi(\bsx)/p(\bsx) \le\kappa<\infty$.
%
\begin{lemma}\label{lem:slicecouple}
%Let $\phi$ be the iterator for the
%Gibbs sampler for distribution $\pi$ on $\xset$.
Let $\pi$ be a density on a bounded rectangular region
$\xset=[\bsa,\bsb]\subset\Real^s$.
Assume that $0<\eta\le\pi(\bsx)\le\kappa<\infty$
holds for all $\bsx\in\xset$.
Let $\xset'=\{(y,\bsx)\mid0\le y\le\pi(\bsx)\}
\subset[\bsa,\bsb]\times[0,\kappa]$
be the domain of the inversive slice sampler.
%with typical state $(\bsx,y)$.
Let $(y_i,\bsx_i)=\phi((y_{i-1},\bsx_{i-1}),\bsu_i)$ for $\bsu
_i\in[0,1]^{s+1}$
be the update for the inversive slice sampler
and put $(y_i',\bsx_i')=\phi((y_{i-1}',\bsx_{i-1}'),\bsu_i)$.
If $\bsu_i\in\Ccal= [0,\eta/\kappa]\times[0,1]^s$,
then $\bsx_i=\bsx_i'$.
\end{lemma}
\begin{pf}
If $u_{i,1}\le\eta/\kappa$,
%in $\Hcal$.
then $y_i=u_{i1}\pi(\bsx_{i-1})$
and $y_i'=u_{i1}\pi(\bsx_{i-1}')$
are in the set $[0,\eta/\kappa]$.
The distribution of $\bsx$ given $y$
for any $y\in[0,\eta/\kappa]$ is $U[\bsa,\bsb]$.
Therefore, $\bsx_i=\bsx_i'
=\bsa+ u_{2\dvtx(s+1)}(\bsb-\bsa)$ (componentwise).
\end{pf}

Lemma \ref{lem:slicecouple} does not couple
the chains because $y_i$ and $y_i'$ are
different in general. But because $\bsx_i=\bsx_i'$,
a coupling will happen at the next step,
that is, $(y_{i+1},\bsx_{i+1})=(y'_{i+1},\bsx'_{i+1})$
when $\bsu_i\in[0,\eta/\kappa]\times[0,1]^s$.
One could revise Theorem \ref{thm:withhomecouple} to
include couplings that happen within some number
$t$ of steps after $\bsu\in\Ccal$ happens. In this
case, it is simpler to say that the chain
whose update comprises two iterations of the inversive
slice sampler satisfies Theorem \ref{thm:withhomecouple}.
For a chain whose update is just one iteration,
the averages over odd and even numbered iterations
both converge properly and so that chain is also consistent.
Alternatively, we could modify the space of $y$
values so that all $y\in[0,\eta/\kappa]$ are identified
as one point. Then $\Ccal$ is a coupling region.

The result of Lemma \ref{lem:slicecouple} also
applies to slice samplers that sample $y\mid\bsx$
and then $\bsx\mid y\sim U\{\bsx\mid\pi(\bsx)\le y\}$
using an $s$-dimensional generator that is not necessarily
inversion.

%CUD sampling of the simple slice sampler
%is consistent when $\xset$ is a Cartesian
%product and
%$0<\varepsilon\le\pi(\bsx)\le M<\infty$
%holds for all $\bsx\in\xset$.

%maybe reversible jump MCMC provides another example

%s5 ###
\section{Riemann integrability}\label{sec:riemaninteg}

Theorem \ref{thm:withhomecouple} proves that
MCMC consistently samples $\pi$ when
implemented using CUD sequences.
We required certain Riemann integrability conditions in defining
regular Rosenblatt--Chentsov transformations.
Here, we verify that nontrivial MCMC algorithms
can have regular Rosenblatt--Chentsov transformations.

It seems odd to use the Riemann integral over
$100$ years after Lebesgue~\cite{lebe1902}.
But pseudo-random number generators
are now typically designed to meet an equidistribution
criterion over rectangular regions \cite{matsnish1998}.
Other times they are designed with a spectral
condition in mind. This again is closely
related to Riemann integrability via the
Weyl \cite{weyl1916} condition where
$\hat\theta_n(f)\to\theta(f)$ for all trigonometric
polynomials $f(\bsx)=e^{2\pi\sqrt{-1}k'\bsx}$
if and only if $\bsx_1,\ldots,\bsx_n$ are
uniformly distributed.
Unless one is using physical random numbers,
the Riemann integral, or perhaps
the improper Riemann integral is
almost implicit.

%s5.1 ###
\subsection{Definitions and basic theorems}

A function from $A\subset\Real^d$ to $\Real^s$ for \mbox{$s\ge1$}
is Riemann integrable if all of its $s$ components are.
To study how Riemann integrability propagates, we will use
the following two definitions.
\begin{definition}
For a function $f\dvtx\Real^k\to\Real$,
the discontinuity set of $f$ is
\[
D(f) = \{\bsx\in\Real^k\mid
f\mbox{ discontinuous at $\bsx$}\}.
\]
If $f$ is only defined on $A\subset\Real^k$,
then $D(f)=D(f_0)$ where
$f_0(\bsx)=f(\bsx)$ for $\bsx\in A$
and $f_0(\bsx)=0$ for $\bsx\notin A$.
\end{definition}
\begin{definition}
For a function $f\dvtx\Real^k\to\Real$, the graph of $f$ is
\[
G(f)=\{(\bsx,y)\in\Real^{k+1}\mid y=f(\bsx)\}.
\]
\end{definition}

Lebesgue's theorem, next, provides a checkable
characterization of Riemann integrability.
\begin{theorem}[(Lebesgue's theorem)]\label{thm:lebesgues}
Let $A\subset\Real^d$ be bounded and let $f\dvtx A\to\Real$
be a bounded function.
%Extend $f$ to all of $\Real^d$ by
%letting it be zero at points not contained in $A$.
Then $f$ is Riemann integrable iff
$\lambda_d(D(f))=0$.
%the points at which
%the extended $f$ is discontinuous form a set of measure
%$0$.
\end{theorem}
\begin{pf}
See Marsden and Hoffman \cite{marshoff1993}, page 455.
\end{pf}

%s5.2 ###
\subsection{Need for Riemann integrable proposals}

Here, we show that Riemann integrability
adds a special requirement to the way an algorithm is
implemented. Then we give an example to show that propagation
rules for Riemann integrability are more complicated than
are those for continuity and differentiability.

Suppose that $F$ is the
$\Ncal\bigl(
{0\choose0},
\bigl({1\atop\rho}\enskip{\rho\atop1}\bigr)
\bigr)$ distribution for some $\rho\in(-1,1)$.
If we take
\[
x_1(\bsu) = \Phi^{-1}(u_1)
\]
and
\[
x_2(\bsu) = \rho x_1(\bsu) + \sqrt{1-\rho^2} \Phi^{-1}(u_2),
\]
then we find that
$f(\bsu) =1_{ a_1\le x_1(\bsu)\le b_1}\times1_{ a_2\le x_2(\bsu
)\le b_2}$
%= 1_{\Phi(z_1-\delta)\le u_1 \le\Phi(z_1+\delta)}
%}.
%$$
%This function
is discontinuous only on a set of measure zero.
It is trivially bounded, and these
two facts imply it is Riemann integrable on $[0,1]^2$.

Another %Rosenblatt
transformation for the same distribution $F$ is
\[
x_1  = \Phi^{-1}(u_1)
\]
and
\[
x_2 =
\cases{
\rho x_1(\bsu) + \sqrt{1-\rho^2} \Phi^{-1}(u_2),&\quad $u_1\notin
\mathbb{Q}$,\cr
-\rho x_1(\bsu) - \sqrt{1-\rho^2} \Phi^{-1}(u_2),&\quad $u_1\in
\mathbb{Q}$.}
\]
Changing the conditional distribution of $x_2$ given $x_1$
on a set of measure $0$ leaves the distribution $F$ of $\bsx$
unchanged. But for this version,
%it yields a Rosenblatt transformation
%that is not regular and
we find
$f$ can be discontinuous on more than a set of measure $0$
and so this inverse Rosenblatt transformation of $F$ is not regular.

In practice, of course, one would use the regular version of the
transformation.
But propagating Riemann integrability to a function built
up from several other functions is not always
straightforward.
%
%A continuous function of a continuous function gives
%a continuous composite function. The same happens for
%differentiable functions.
The core of the problem is that the
composition of two Riemann integrable functions
need not be Riemann integrable.

As an example \cite{ghorlima2006},
consider Thomae's function on $(0,1)$,
\[
f(x) =
\cases{
1/q, &\quad $x= p/q\in\Rational$,\cr
0, &\quad else,}
\]
where it is assumed that $p$ and $q$ in the representation $p/q$
have no common factors.
This $f$ is continuous except on $\Rational\cap(0,1)$ and so
it is Riemann integrable. The function $g(x)=1_{0<x\le1}$ is
also Riemann integrable. But $g(f(x))=1_{x\in\Rational}$
for $x\in(0,1)$,
which is famously not Riemann integrable.
The class of
Riemann integrable functions, while more
restrictive than we might like for conclusions,
is also too broad to use in propagation rules.

%s5.3 ###
\subsection{Specializing to MCMC}

First, we show that the acceptance-rejection step
in Metropolis--Hastings does not cause problems
with Riemann integrability.
\begin{lemma}\label{lem:accrej}
Let $k\in\Natural$ and suppose
that $g$, $h$ and $A$ are real-valued Riemann integrable
functions on $[0,1]^k$.
For $\bsu\in[0,1]^{k+1}$ define
\[
f(\bsu) = \cases{
g(\bsu_{1\dvtx k}),&\quad $u_{k+1}\le A(\bsu_{1\dvtx k})$,\cr
h(\bsu_{1\dvtx k}),&\quad else.}
\]
Then $f$ is Riemann integrable on $[0,1]^{k+1}$.
\end{lemma}
\begin{pf}
First,
$D(f)\subset((D(g)\cup D(h))\times[0,1])\cup G(A)$.
Riemann integrability of $g$ gives $\lambda_k(D(g))=0$.
Similarly, $\lambda_k(D(h))=0$.
Therefore, $\lambda_{k+1}(D(g)\cup D(h))\times[0,1])=0$.

Turning to $G(A)$, we % may
split the domain $[0,1]^k$ of $A$ into
$n^k$ congruent subcu\-bes $C_{n,1},\ldots,C_{n,n^k}$
(whose boundaries overlap).
Then $G(A)\subseteq\bigcup_{i=1}^{n^k} C_{n,i}\times[m_{i,n},M_{i,n}]$,
where $m_{i,n} = \inf_{\bsu_{1\dvtx k}\in C_{n,i}}A(\bsu_{1\dvtx k})$
and $M_{i,n}=\sup_{\bsu_{1\dvtx k}\in C_{n,i}}A(\bsu_{1\dvtx k})$.
As a result $\lambda_{k+1}(G(h))\le n^{-k}\sum_i(M_{i,n}-m_{i,n})$.
Riemann integrability
of $A$ implies this upper bound vanishes as $n\to\infty$.
Therefore, $\lambda_{k+1}(G(A)) = 0$ and so $\lambda_{k+1}(D(f))=0$
and the result follows by Lebesgue's theorem.
\end{pf}

In the MCMC context, $g$ and $h$ are the $j$th component
of the proposal and the previous state, respectively,
$A$ is the acceptance probability, and $\bsu$ is the
ensemble of uniform random variables used in $m$ stage
Rosenblatt--Chentsov coupling and $k=(m+1)d-1$.

For consistency results, we
study the proportion of times $f(\bsu)\in[\bsa,\bsb]$.
It is enough to consider the components
one at a time and in turn to show
$1_{f_j(\bsu)\le b_j}$ and
$1_{f_j(\bsu)< a_j}$ are Riemann integrable.
However, as the example with Thomae's function shows,
even the indicator function of an interval
applied to a Riemann integrable function can
give a non-Riemann integrable composite function.

We may avoid truncation by employing
bounded continuous test functions.
We will use the following simple corollary
of Lebesgue's theorem.
\begin{lemma}\label{lem:ctsofriemann}
For $k\ge1$ and $r\ge1$, let $g_1,\ldots,g_r$ be Riemann integrable functions
from $[0,1]^k$ to a bounded interval $[a,b]\subset\Real$.
Let $h$ be a continuous function from $[a,b]^k$ to $\Real$.
Then
\[
f(\bsu)=h(g_1(\bsu),\ldots,g_k(\bsu))
\]
is Riemann integrable on $[0,1]^k$.
\end{lemma}
\begin{pf}
Because $h$ is continuous,
$D(f)\subset\bigcup_{j=1}^r D(g_k)$.
But \mbox{$\lambda_k(D(g_k)) =0$}. Therefore, $\lambda_k(D(f))=0$
and so $f$ is Riemann integrable by Lebesgue's \mbox{theorem}.
\end{pf}

We can also propagate Riemann
integrability through monotonicity.
If $g$ is a monotone function from $\Real$ to $\Real$
and $f$ is the indicator of an interval, then
$f\circ g$ is the indicator of an interval too,
and hence is Riemann integrable, when
that interval is of finite length.
\begin{lemma}\label{lem:rosenblatt}
Let $F_1(x_1)$ be the CDF of $x_1$
and for $j=2,\ldots,s$, let $F_j(x_j\mid\bsx_{1\dvtx(j-1)})$ be the
conditional CDF of $x_j$ given $\bsx_{1\dvtx(j-1)}$.
%$F_j(u_j\mid x_1,\ldots,x_{j-1})$
%is piecewise in $(x_1,\ldots,x_{j-1})$
Suppose that the CDFs $F_j(x_j\mid\bsx_{1\dvtx(j-1)})$
are continuous functions of $\bsx_{1\dvtx j}$ and
that the quantile functions
$F_j^{-1}(u_j\mid\bsx_{1\dvtx (j-1)})$
are continuous in $(u_j,\bsx_{1\dvtx(j-1)})\in[0,1]\times\Real^{j-1}$,
for $j=2,\ldots,s$.
Define functions
$z_1(\bsu)=F_1^{-1}(u_1)$ and
$z_j(\bsu)=F_j^{-1}(u_j\mid\mathbf{z}_{1\dvtx(j-1)}(\bsu))$
for $j=2,\ldots,s$, where $\mathbf{z}_{1\dvtx(j-1)} = (z_1,\ldots, z_{j-1})$.
Then for \mbox{$\bsb\in\Real^s$}, the set
\[
S(\bsb)=\{\bsu\mid z_j(\bsu) \le b_j, 1\le j\le s\}
\]
is Jordan measurable.
\end{lemma}
\begin{pf}
%For $k=1,\ldots,s$, let
%$S_k = \{ \bsu_{1:k}\in[0,1]^k \mid x_j\le b_j, 1\le j\le k\}$.
%Define functions
%$x_1(\bsu)=F_1^{-1}(u_1)$ and
%$x_j(\bsu)=F_j^{-1}(u_j\mid\bsx_{1:(j-1)}(\bsu))$
%for $j=2,\ldots,s$.
By hypothesis, $z_k$ is a continuous
function of $\bsu\in[0,1]^s$, for $k=1,\ldots,s$, and so is
$F_k(b_k\mid\mathbf{z}_{1\dvtx(k-1)}(\bsu))$.
This latter only depends on $\bsu_{1\dvtx(k-1)}$,
for $k=2,\ldots,s$, and so we write it as
$g_k(\bsu_{1\dvtx(k-1)})$.

For $k=1,\ldots,s$, let $S_k = \{\bsu_{1\dvtx k}\mid u_j \le g_j(\bsu
_{1\dvtx (j-1)}) \mbox{ for } j = 1,\ldots, k\}$.
The set $S_1$ is the interval $[0,F_1^{-1}(b_1)]$,
and hence is Jordan measurable.
Suppose $S_k$ is Jordan measurable for $k<s$.
Then
\begin{eqnarray*}
S_{k+1} & =
( S_k\times[0,1]) \cap G_{k+1} \qquad
\mbox{where }
G_{k+1} =\bigl\{ \bsu_{1\dvtx (k+1)}\mid u_{k+1} \le g_{k+1} (\bsu
_{1\dvtx k})\bigr\}.
\end{eqnarray*}
The set $S_k\times[0,1]$ is Jordan measurable because $S_k$ is.
The boundary of $G_{k+1}$ is contained within the
intersection of the graph of $g_{k+1}$ and the boundary
of $[0,1]^{k+1}$ and so $G_{k+1}$ is Jordan measurable.
The result follows by induction because $S(\bsb) = S_s$.
\end{pf}

%s5.4 ###
\subsection{Regularity of Rosenblatt--Chentsov}

Here, we give sufficient conditions for the Rosenblatt--Chentsov
transformation to be regular.
\begin{theorem}
For integer $m\ge0$, let $\bsx_m$
be the endpoint of the Rosenblatt--Chentsov transformation
of $[0,1]^{(d+1)m}$, started with
a Riemann integrable function $\psi_R$
and continued via the Metropolis--Hastings
update $\phi$.
Let $\phi$ be defined in terms of the
proposal function $\bsy\dvtx \Real^s\times[0,1]^{d-1}\to\Real^s$
with proposal density $p(\cdot,\cdot)\dvtx \Real^s\times\Real^s\to
[0,\infty)$
and target density $\pi\dvtx \Real^s\to[0,\infty)$.
Let $f$ be a bounded continuous function
%with bounded support set in
on $\Real^s$.

If $\psi$ is bounded and $\bsy$, $P$
and $\pi$ are bounded continuous
functions, then
$f(\bsx_m(\bsu_0,\ldots,\bsu_m))$ is a Riemann integrable function
of the variables $[0,1]^{(d+1)m}$ used in
the Rosenblatt--Chentsov transformation.
\end{theorem}
\begin{pf}
We only need to show that
$\bsx_{m}$ is a Riemann integrable function
of $(\bsu_0,\ldots,\bsu_m)\in[0,1]^{d(m+1)}$
and then the result follows by
Lemma \ref{lem:ctsofriemann}.

We proceed by induction.
For $m=0$, $\bsx_0=\psi(\bsu_0)$
is bounded and continuous on $[0,1]^d$, hence
it is Riemann integrable.

Now suppose that $\bsx_{m-1}$ is a Riemann
integrable function on $[0,1]^{dm}$.
Let $h(\bsu_0,\ldots,\bsu_{m-1},\bsu_{m 1\dvtx (d-1)})$
be the value $\bsx_{m-1}$, written as a Riemann
integrable function on $[0,1]^{dm+d-1}$, so it
ignores its last $d-1$ arguments.
Let $g(\bsu_0, \ldots,\bsu_{m-1},\bsu_{m 1\dvtx (d-1)})$
be the proposal
$\bsy_m=\bsy(\bsx_{m-1},\bsu_{m 1\dvtx (d-1)})
=\bsy(g(\cdot),\bsx_{m-1}$, $\bsu_{m 1\dvtx (d-1)})$.
This is a continuous function $\bsy(\cdot,\cdot)$ of two Riemann integrable
functions on $[0,1]^{d(m+1)-1}$ and so it is Riemann
integrable.
Next, $A(\cdot,\cdot)$ is a continuous function of
both $\bsx_{m-1}$ and $\bsy_m$ which are in turn Riemann
integrable functions on $[0,1]^{dm+d-1}$, and so $A(\cdot,\cdot)$
is Riemann integrable.
Then $\bsx_m$ is a Riemann integrable function
on $[0,1]^{dm+d}$, by Lemma~\ref{lem:accrej}, completing
the induction.
\end{pf}

%s6 ###
\section{Conditions for the Gibbs sampler}\label{sec:gibbsexample}

In studying the Gibbs sampler, we made several assumptions.
First, we required Jordan measurability for the
sets $\Bcal_m(\bsx,\widehat\bsx)$. Second, we required
a contraction property. In this section, we show
that those assumptions are reasonable.

%s6.1 ###
\subsection{Jordan measurability of $\Bcal_m(\bsx,\widehat\bsx)$}

We give an example where the conditions of
Theorem \ref{thm:withiterfunction} are satisfied, that is, the sets
$\Bcal_{m}(\bsx,\widehat{\bsx})$ are Jordan measurable for all $m
\ge1$ and $\bsx,\widehat{\bsx} \in\Omega$ (for some suitable
domain $\Omega\subset\Real^s$). Assume (additionally to the
assumptions made in Theorem \ref{thm:withiterfunction}) that
$\phi(\bsx,\bsu)$ is totally differentiable with continuous
derivative with respect to $\bsu$ for each $\bsx\in\Omega$ and
that $d$ is based on the $L_p$ norm for some $1\le p < \infty$.
Further, assume that the gradient of
$d(\phi(\bsx,\bsu),\phi(\widehat{\bsx},\bsu))$ with respect to
$\bsu$ vanishes only on a null set for all $\bsx, \widehat{\bsx}
\in
\Omega$, $\bsx\neq\widehat{\bsx}$, that is,
\[
\lambda\bigl(\{\bsu\in[0,1]^{d}\dvtx  \nabla_{\bsu}\, d(\phi(\bsx,\bsu),\phi(\widehat{\bsx},\bsu)) = \bszero
\}\bigr) = 0,
\]
for all $\bsx, \widehat{\bsx} \in\Omega$, $\bsx\neq
\widehat{\bsx}$, where $\lambda$ denotes the Lebesgue measure and
where $\nabla_{\bsu}\, d(\phi(\bsx,\bsu),\phi(\widehat{\bsx},\bsu
)) =
(\frac{\partial}{\partial u_j}\,
d(\phi(\bsx,\bsu),\phi(\widehat{\bsx},\bsu)))_{j=1,\ldots
, d}$
denotes the gradient.

Then, for all $m \ge1$, we also have
\[
\lambda\bigl(\{\bsu\in[0,1]^{dm}\dvtx  \nabla_{\bsu}\,
d(\phi_m(\bsx,\bsu),\phi_m(\widehat{\bsx},\bsu)) = \bszero \}\bigr) = 0
\]
for all $\bsx,\widehat{\bsx} \in\Omega$, $\bsx\neq\widehat{\bsx
}$. Let $\bsx, \widehat{\bsx}
\in\Omega$ with $\bsx\neq\widehat{\bsx}$ be fixed. Then for
almost all $\bsu^\ast\in[0,1]^{dm}$
we have $\nabla_{\bsu}\, d(\phi_m(\bsx,\bsu^\ast),\phi_m(\widehat{\bsx},\bsu^\ast))
\neq\bszero$. Therefore, there is a $\delta> 0$ such that $\nabla
_{\bsu}
d(\phi_m(\bsx,\bsu),\phi_m(\widehat{\bsx},\bsu)) \neq\bszero$ for
all $\bsu\in N_\delta(\bsu^\ast)$, where $N_{\delta}(\bsu^\ast) =
\{\bsv\in[0,1]^{dm}\dvtx  \|\bsu^\ast-\bsv\|_{L_2} < \delta\}$ is a
neighborhood of $\bsu^\ast$. Therefore, the directional derivative at
a point $\bsu\in N_\delta(\bsu^\ast)$ is different from $0$, except
on a hyperplane, that is, almost everywhere. Hence, by the mean value
theorem, the function
$d(\phi_m(\bsx,\bsu),\phi_m(\widehat{\bsx},\bsu))$ for $\bsu\in
N_\delta(\bsu^\ast)$ can at most be constant on a hyperplane, which
has Lebesgue measure $0$. Note that $N_\delta(\bsu^\ast) \cap
\mathbb{Q}^{dm} \neq\varnothing$, therefore there is a countable
number of elements $\bsu^\ast_1, \bsu^\ast_2, \ldots$ and numbers
$\delta_1,\delta_2,\ldots$ with the properties of $\bsu^\ast$ and
$\delta$ described above and for which we have $\bigcup_{n=1}^\infty
N_{\delta_n}(\bsu^\ast_n) = [0,1]^{dm}$. Therefore, we have
\[
\lambda\bigl(\{\bsu\in[0,1]^{dm}\dvtx
d(\phi_m(\bsx,\bsu),\phi_m(\widehat{\bsx},\bsu)) = c \}
\bigr)
= 0,
\]
for any $c > 0$.

The set of points where $1_{\Bcal_{m}(\bsx,\widehat{\bsx})}$ is
discontinuous is given by
\begin{eqnarray*}
D & = & \{\bsu\in[0,1]^{dm}\dvtx  \forall\delta> 0\ \exists\bsv,\bsv'
\in N_\delta(\bsu) \mbox{ such that } \\
&&\hspace*{4.5pt}
d(\phi_m(\bsx,\bsv),\phi_m(\widehat{\bsx},\bsv)) > \gamma^m
\mbox{
and } d(\phi_m(\bsx,\bsv'),\phi_m(\widehat{\bsx},\bsv')) \le
\gamma^m \}.
\end{eqnarray*}
As $\Bcal_{m}(\bsx,\widehat{\bsx})$ and $\{\bsu\in[0,1]^{dm}\dvtx
d(\phi_m(\bsx,\bsu),\phi_m(\widehat{\bsx},\bsu)) < \gamma^m
\}$ are open, it follows that
\[
D \subseteq\{\bsu\in[0,1]^{dm}\dvtx
d(\phi_m(\bsx,\bsu),\phi_m(\widehat{\bsx},\bsu)) = \gamma^m
\}.
\]
Therefore, $\lambda_{dm}(D) = 0$ and Lebesgue's theorem (see
Theorem \ref{thm:lebesgues}) implies that
$\Bcal_{m}(\bsx,\widehat{\bsx})$ is Jordan measurable.

%s6.2 ###
\subsection{Contraction}

Here, we illustrate how the
Gibbs sampler yields a~contraction for the probit model.
In this model,
\[
Z_i  = \bsx_i^\tran\beta+ \epsilon_i
\]
and
\[
Y_i  = \bsone_{Z_i>0},
\]
for $i=1,\ldots,n$ for independent
$\epsilon_i \sim\Ncal(0,1)$. The coefficient $\beta\in\Real^p$
has a~noninformative prior distribution.
The predictors are $\bsx_i\in\Real^p$. We define the matrix $X$ with
$ij$ element $x_{ij}$.
We assume that $X$ has rank $p$.

The state of the Markov chain
is $(\beta,\bsZ)\in\Omega\subset\Real^{p+n}$,
where $\bsZ= (Z_1,\ldots,\break Z_n)^\tran$.
Given the observed data $(y_1,\ldots,y_n,\bsx_1,\ldots,\bsx_n)$, we
can use the Gibbs sampler to simulate
the posterior distribution of $\beta$ and $\bsZ= (Z_1,\ldots,Z_n)^\tran$.
A single step of the Gibbs sampler makes the transition
\[
\pmatrix{\beta^{(k-1)}\cr\bsZ^{(k-1)}}
\stackrel{u_1,\ldots,u_n}{\longrightarrow}
\pmatrix{\beta^{(k-1)}\cr\bsZ^{(k)}}
\stackrel{u_{n+1},\ldots,u_{n+p}}{\longrightarrow}
\pmatrix{\beta^{(k)}\cr\bsZ^{(k)}}
\]
for $k\ge1$ using generators given explicitly below.
The values $u_1,\ldots,u_{n+p}$ are the components
of $\bsu_k\in(0,1)^{n+p}$.
We also write the transitions as
\[
(\beta,\bsZ)\to
\phi( (\beta,\bsZ),\bsu)
=\bigl(\phi^{(1)}( (\beta,\bsZ),\bsu),\phi^{(2)}( (\beta,\bsZ
),\bsu)\bigr),
\]
where $\phi$ and its components $\phi^{(1)}$ and $\phi^{(2)}$ (for
$\beta$
and $\bsZ$, resp.) are given explicitly below.

Given $\beta$, the components of $\bsZ$ are independent, with
\[
Z_i \sim
\cases{
\Ncal(\bsx_i^\tran\beta,1)|Z_i>0, &\quad if $Y_i = 1$,\vspace*{2pt}\cr
\Ncal(\bsx_i^\tran\beta,1)|Z_i\leq0, &\quad if $Y_i = 0$.}
\]
We may generate them from $u_1,\ldots,u_n\in(0,1)$ by
%
%e9 ###
%
\begin{equation}\label{eq:zfrombeta}
Z_i =
\cases{
\bsx_i^\tran\beta+\Phi^{-1}\bigl(\Phi(-\bsx_i^\tran\beta
)+u_i\Phi(\bsx_i^\tran\beta)\bigr), &\quad if
$Y_i=1$,\vspace*{2pt}\cr
\bsx_i^\tran\beta+\Phi^{-1}(u_i\Phi(-\bsx_i^\tran\beta
)), &\quad if
$Y_i=0$.}
\end{equation}

Given $\bsZ$, the distribution of $\beta$ is
$\beta\sim\Ncal((X^\tran X)^{-1}X^\tran\bsZ,(X^\tran X)^{-1})$.
We may generate it using $u_{n+1},\ldots,u_{n+p}\in(0,1)$
via
%
%e10 ###
%
\begin{equation}\label{eq:betafromz}
\beta= (X^\tran X)^{-1}X^\tran\bsZ+ (X^\tran X)^{-1/2}
\pmatrix{\Phi^{-1}(u_{n+1})\cr\vdots\cr\Phi^{-1}(u_{n+p})}.
\end{equation}
Thus equation (\ref{eq:betafromz}) defines $\phi^{(1)}$ while
(\ref{eq:zfrombeta}) defines $\phi^{(2)}$.

The framework in \cite{alsmfuh2001} allows one to pick
a metric that conforms to the problem.
We use the metric
$
d((\beta,\bsZ),(\beta',\bsZ'))
=\max(d_1(\beta,\beta'),d_2(\bsZ,\bsZ')),
$
where
%
%e11 ###
%
\begin{equation}
\label{eq:d1}
d_1(\beta,\beta') = d_1(\beta-\beta')=\sqrt{(\beta-\beta
')^\tran(X^\tran X)(\beta- \beta')}
\end{equation}
and
%
%e12 ###
\begin{equation}
\label{eq:d2}
d_2(\bsZ,\bsZ') = d_2(\bsZ-\bsZ')=\sqrt{(\bsZ-\bsZ')^\tran
(\bsZ- \bsZ')}.
\end{equation}
We show below that
%
%e13 ###
%
\begin{equation}\label{eq:itcontracts}\quad
d\bigl(\bigl(\beta^{(k)},\bsZ^{(k)}\bigr),\bigl(\beta^{\prime(k)},\bsZ^{\prime
(k)}\bigr)\bigr)
\le d\bigl(\bigl(\beta^{(k-1)},\bsZ^{(k-1)}\bigr),\bigl(\beta^{\prime
(k-1)},\bsZ^{\prime(k-1)}\bigr)\bigr)
\end{equation}
for pairs
$(\beta^{(k-1)},\bsZ^{(k-1)}),(\beta^{\prime(k-1)},\bsZ^{\prime
(k-1)})$
of distinct points in $\Omega$.
%where $\phi$ is the transition function of the Gibbs sampler and $\bsu
Both metrics $d_1$ and $d_2$ are also norms, which simplifies our task.
%We also write the output $\phi((\beta,\bsZ),\bsu)$
%as $( \phi_\beta((\beta,\bsZ),\bsu),\phi_{\bsZ}((\beta,\bsZ),\bsu))$
%where $\phi_\beta(\cdot)$ and $\phi_{\bsZ}(\cdot)$ extract their
%named components.

Suppose first that $\beta^{(k-1)}=\beta^{\prime(k-1)}$.
Then it follows easily that $\bsZ^{(k)}=\bsZ^{\prime(k)}$
and $\beta^{(k)}=\beta^{\prime(k)}$,
so then the left-hand side of (\ref{eq:itcontracts}) is $0$.
%$\phi_{\bsZ}((\beta,\bsZ),\bsu)=\phi_{\bsZ}((\beta',\bsZ'),\bsu)$
%and then
%$\phi_{\beta}((\beta,\bsZ),\bsu)=\phi_{\beta}((\beta',\bsZ'),\bsu)$
%so that
%$d(\phi((\beta,\bsZ),\bsu),\phi((\beta',\bsZ'),\bsu))=0
%<d((\beta,\bsZ),(\beta',\bsZ'))=0$
As a result, we may assume without loss of generality that
$d_1(\beta^{(k-1)}-\beta^{\prime(k-1)})>0$.
With this assumption, we will use the bound
%
%e14 ###
%
\begin{eqnarray}
&&\frac{d((\beta^{(k)},\bsZ^{(k)}),(\beta^{\prime(k)},\bsZ
^{\prime(k)}))}
{d((\beta^{(k-1)},\bsZ^{(k-1)}),(\beta^{\prime(k-1)},\bsZ
^{\prime(k-1)}))}\nonumber\\[-8pt]\\[-8pt]
&&\qquad\le\max\biggl(
\frac{d_1(\beta^{(k)}-\beta^{\prime(k)})}
{d_1(\beta^{(k-1)}-\beta^{\prime(k-1)})},
\frac{d_2(\bsZ^{(k)}-\bsZ^{\prime(k)})}
{d_1(\beta^{(k-1)}-\beta^{\prime(k-1)})}
\biggr).\nonumber
\end{eqnarray}

We begin by studying the update to $\bsZ$.
Subtracting $\bsx_i^\tran\beta$
from both sides of (\ref{eq:zfrombeta}), applying $\Phi(\cdot)$,
differentiating with
respect to $\beta$ and gathering up terms, we find that
$\frac{\partial}{\partial\beta} Z_i = \lambda_i\bsx_i$
where
%
%e15 ###
%
\begin{equation}\label{eq:lambdafrom}
\lambda_i =
\cases{
1 - \dfrac{(1-u_i)\varphi(\bsx_i^\tran\beta)}{\varphi
(Z_i-\bsx_i^\tran\beta)}, &\quad if
$Y_i=1$,\vspace*{2pt}\cr
1 - \dfrac{u_i\varphi(-\bsx_i^\tran\beta)}{\varphi
(Z_i-\bsx_i^\tran\beta)}, &\quad if
$Y_i=0$,}
\end{equation}
and $\varphi$ is the $\Ncal(0,1)$ probability density function.

It is clear that $\lambda_i < 1$. Next, we show that $\lambda_i\ge0$.
We begin by inverting~(\ref{eq:zfrombeta}) to get
%
%e16 ###
%
\begin{equation}\label{eq:ufrombeta}
u_i =
\cases{
\dfrac{\Phi(Z_i-\bsx_i^\tran\beta)-\Phi(-\bsx_i^\tran\beta
)}{\Phi(\bsx_i^\tran\beta)},
&\quad if $Y_i=1$,\cr
\dfrac{\Phi(Z_i-\bsx_i^\tran\beta)}{\Phi(-\bsx_i^\tran\beta)},
&\quad if $Y_i=0$.}
\end{equation}
Substituting (\ref{eq:ufrombeta}) into (\ref{eq:lambdafrom}) and simplifying
yields
%
%e17 ###
%
\begin{equation}\label{eq:lambdafrom2}
1-\lambda_i =
\cases{
\dfrac{\varphi(\bsx_i^\tran\beta)\Phi(-Z_i+\bsx
_i^\tran\beta)}{\Phi(\bsx_i^\tran\beta)\varphi(Z_i-\bsx_i^\tran
\beta)},
&\quad if $Y_i=1$,\vspace*{2pt}\cr
\dfrac{\varphi(-\bsx_i^\tran\beta)\Phi(Z_i-\bsx
_i^\tran\beta)}{\Phi(-\bsx_i^\tran\beta)\varphi(Z_i-\bsx
_i^\tran\beta)},
&\quad if $Y_i=0$.}
\end{equation}
Now consider the function $\tau(x)=\varphi(x)/\Phi(x)$. This
function is nonnegative
and decreasing, using a Mill's ratio bound from \cite{gord1941}.
When $Y_i=1$, then $1-\lambda_i=\tau(\bsx_i^\tran\beta)/\tau(\bsx
_i^\tran\beta-Z_i)\le1$ because
then $Z_i\ge0$. We also used symmetry of $\varphi(\cdot)$.
If instead $Y_i=0$, then $1-\lambda_i=\tau(-\bsx_i^\tran\beta
)/\tau(-\bsx_i^\tran\beta+Z_i)\le1$
because then $Z_i\le0$.
Either way, $1-\lambda_i\le1$ and therefore $\lambda_i\in[0,1)$ for all $i$.

Writing the previous results in a compact matrix form, we have
% \begin{align*}
% \frac{\partial\bsZ}{\partial\beta} & = (\begin{array}{cccc}
% \ds\frac{\partial z_1}{\partial\beta_1} & \ds\frac{\partial z_1}{
% \ds\frac{\partial z_2}{\partial\beta_1} & \ds\frac{\partial z_2}{
% \vdots& \vdots& \ddots& \vdots\\
% \ds\frac{\partial z_n}{\partial\beta_1} & \ds\frac{\partial z_n}{
% \end{array})
% = \Lambda X,
% \end{align*}
%
\[
\frac{\partial\bsZ}{\partial\beta} = \biggl( \ds\frac{\partial
z_i}{\partial\beta_j}\biggr)_{ij} = \Lambda X,
\]
where
%$$X = (\begin{array}{c}
$\Lambda= \Lambda(\beta,\bsZ)=\diag(\lambda_1,\ldots,\lambda_n)$.
Similarly, equation (\ref{eq:betafromz}) yields
\[
\frac{\partial\beta}{\partial\bsZ}=(X^\tran X)^{-1}X^\tran.
\]

%The second step is to update $\beta$ when fixing $\bsZ$. Notice the
%conditional distribution of $\beta$ given $\bsZ$ is multivariate
%normal:
%$$\beta\sim\Ncal((X^\tran X)^{-1}X^\tran\bsZ,(X^\tran X)^{-1})$$
%So we could use p random numbers $(u_{n+1},u_{n+2}...u_{n+p})$ to
%simulate a new $\beta$ from this conditional distribution as follows:
%$$\beta= (X^\tran X)^{-1}X^\tran\bsZ+ f(u_{n+1},u_{n+2}...u_{n+p})$$
%where $f(u_{n+1},u_{n+2}...u_{n+p})$ produces a random vector from
%p-dimensional Gaussion distribution $\Ncal(0,(X^\tran X)^{-1})$.
%Thus it's easy to see that $$\ds\frac{\partial\beta}{\partial\bsZ}=
%(X^\tran X)^{-1}X^\tran$$

Thus, for the $\bsZ$ update with any $\bsu_k\in(0,1)^{n+p}$,
%
%e18 ###
%
\begin{eqnarray}\label{eq:zup}
\frac{d_2(\bsZ^{(k)}-\bsZ^{\prime(k)})}
{d_1(\beta^{(k-1)}-\beta^{\prime(k-1)})}
&\le&\mathop{\sup_{\wt\beta^{(k-1)},\wt\bsZ^{(k)}}}_{d_1(\xi
) = 1}d_2\biggl(\frac{\partial\wt\bsZ^{(k)}}{\partial\wt\beta
^{(k-1)}} \xi\biggr)\nonumber\\[-8pt]\\[-8pt]
&\le&{\mathop{\sup_{\beta,\bsZ}}_{(X\xi)^\tran X\xi= 1}}
\Vert\Lambda(\beta,\bsZ) X \xi\Vert
<1.\nonumber
\end{eqnarray}

%Summing up, the Gibbs Sampler algorithm has a 2-step procedure which
%is illustrated as follows:
%$$\beta^{(k)}\stackrel{u_1...u_n}{\longrightarrow}\bsZ^{(k)}

For the $\beta$ update, applying the
chain rule gives
\[
\frac{\partial\beta^{(k)}}{\partial\beta^{(k-1)}} = \frac
{\partial\beta^{(k)}}{\partial\bsZ^{(k-1)}}\frac{\partial\bsZ
^{(k-1)}}{\partial\beta^{(k-1)}}
= (X^\tran X)^{-1}X^\tran\Lambda X
\]
and then %(noting that $\beta^{(k)}$ depends on $\bsZ^{(k)}$ but not $
%
%e19 ###
%
\begin{eqnarray}\label{eq:bup}
\frac{d_1(\beta^{(k)}-\beta^{\prime(k)})}
{d_1(\beta^{(k-1)}-\beta^{\prime(k-1)})}
%&\le\sup_{\stackrel{\ds\wt\beta}{d_1(\xi) = 1}}d_1(\frac{
&\le&\sup_{\wt\beta, d_1(\xi) = 1}d_1\biggl(\frac{\partial\wt
\beta^{(k)}}{\partial\wt\beta^{(k-1)}} \xi\biggr)\nonumber\\
&=&\mathop{\sup_{\beta,\bsZ}}_{(X\xi)^\tran X\xi= 1}
d_1((X^\tran X)^{-1}X^\tran\Lambda X\xi)\nonumber\\
&=&\sup_{\beta,\bsZ,\Vert\eta\Vert= 1}
d_1((X^\tran X)^{-1}X^\tran\Lambda\eta)\nonumber\\[-8pt]\\[-8pt]
&=&{\sup_{\beta,\bsZ,\Vert\eta\Vert= 1}}
\Vert X(X^\tran X)^{-1}X^\tran\Lambda\eta\Vert
\nonumber\\
& \le&{\max_{1\le i\le n}\lambda_i} \nonumber\\
& <&1,\nonumber
\end{eqnarray}
using the nonexpansive property
of the projection matrix $X(X^\tran X)^{-1}X^\tran$.
%Plugging the derived equations of $\frac{\partial\beta^{(k+1)}}{
%we get:
%X)^{-1}X^\tran\Lambda X

By combining (\ref{eq:zup}) with (\ref{eq:bup}), we establish
the contraction (\ref{eq:itcontracts}).

\section{Open versus closed intervals}\label{sec:openclose}

In the Lebesgue formulation,
$U(0,1)^d$ and $U[0,1]^d$
are the same distribution, in that they cannot
be distinguished with positive probability
from any countable sample of independent values.
Riemann integrals are usually defined for $[0,1]^d$
and discrepancy measures are usually
defined for either $[0,1]^d$ or $[0,1)^d$.
These latter theories are designed
for bounded functions.

In Monte Carlo simulations, sometimes values
$u_{ij}\in\{0,1\}$ are produced. These end points
can be problematic with inversion, where they may
yield extended real values, and hence good practice
is to select random number generators supported
in the open interval $(0,1)$.

For our Gibbs sampler example with the probit
model, we required $\bsu_k\in(0,1)^{n+p}$.
This was necessary because otherwise the
values $\phi(\bsx,\bsu)$ might fail to belong
to $\xset$.

Our slice sampler example had $\xset$
equal to the bounded rectangle $[\bsa,\bsb]$.
Then values $u_{ij}\in\{0,1\}$ do not generate sample
points outside $\xset$.

Our Metropolized independence sampler did not
require bounded support.
It could produce extended real values.
Those however are not problematic
for weak convergence, which is based
on averages of $1_{[\bsa,\bsb]}(\bsx_i)$
or other bounded test functions.
Also, the chain will not get stuck at an unbounded
point.

%s8 ###
\section{Discussion}\label{sec:conclusions}

We have demonstrated that MCQMC algorithms
formed by Metropolis--Hastings updates driven
by completely uniformly distributed points can
consistently sample a continuous stationary distribution.
Some regularity conditions are required, but
we have also shown that those conditions hold
for many, though by no means all, MCMC updates.
The result is a kind of ergodic theorem for QMC
like the ones in \cite{chen1967}
and \cite{qmcmetro} for finite state spaces.

When RQMC is used in place of QMC to drive an MCMC
simulation, then instead of CUD points, we need
to use weakly CUD points. These satisfy
$\Pr(D_n^{*d}>\epsilon)\to0$ for all $\epsilon>0$
and all $d\in\Natural$.

Our version of MCMC above leaves out some methods in which one or more
components
of $\bsu_i$ are generated by acceptance-rejection sampling
because then we cannot assume $d<\infty$.
A modification based on splicing \IID\  $U[0,1]$ random variables
into a CUD sequence was proposed by Liao \cite{liao1998}
and then shown to result in a weakly CUD sequence in
\cite{qmcmetro2}.

We do not expect that a global substitution of QMC points
will always bring a large improvement to MCMC algorithms.
What we do expect is that means of smooth functions of
the state vector in Gibbs samplers will often benefit
greatly from more even sampling.

It is also a fair question
to ask when one needs an MCMC result computed to
the high precision that QMC sometimes makes possible.
Gelman and Shirley \cite{gelmshir2010} address
this issue, distinguishing Task 1 (inference about
a~parameter $\theta$) from Task 2
[precise determination of $\mathbb{E}(\theta)$
or more generally $\mathbb{E}(f(\theta))$ conditional
on the data, or a posterior quantile of $\theta$].
The accuracy of Task 1 problems may be limited more
by sample size than by Monte Carlo effort.
Task 2 problems include computation of normalizing
constants and problems where one wants to report
numerically stable, and hence more reproducible,
simulation output.

\begin{appendix}\label{app}
\section*{Appendix: Proofs}

This Appendix contains the lengthier proofs.

We need one technical lemma about CUD points.
Consider overlapping blocks of $dk$-tuples from $u_i$,
with starting indices $d$ units apart. If $u_i$
are CUD then these overlapping blocks are uniformly
distributed. The proof works by embedding the
$dk$-tuples into nonoverlapping $rdk$-tuples. For large
$r$, the boundary effect between adjacent blocks
becomes negligible. This result is also needed
for the argument in~\cite{qmcmetro}.
\begin{lemma}\label{lem:tech}
For $j\ge1$, let $u_j\in[0,1]$.
For integers $d,i,k\ge1$, let
$\bsx_i = (u_{d(i-1)+1},\ldots,u_{d(i-1)+dk})$.
If $u_j$ are completely\vspace*{1pt} uniformly distributed,
then $\bsx_i\in[0,1]^{dk}$ are uniformly distributed.
\end{lemma}
\begin{pf}%[\bf Proof of Lemma \ref{lem:tech}.]
Choose any $\bsc\in[0,1]^{dk}$.
Let $v=\prod_{j=1}^{dk}c_j$ be the volume of $[\bszero,\bsc)$.
For integers $r\ge1$, define $f_r$ on $[0,1]^{rdk}$ by
$f_r(\bsu) = \sum_{j=0}^{(r-1)k}
1_{[\bszero,\bsc)}(u_{jd+1},\ldots,\break u_{jd+dk})$.
Each $f_r$ has Riemann integral $((r-1)k+1)v$. We
use $f_r$ on nonoverlapping blocks of length $rdk$ from $u_j$:
\begin{eqnarray*}
&&\frac1n\sum_{i=1}^n1_{[\bszero,\bsc)}(\bsx_i)
\ge\frac1n
\sum_{i=1}^{\lfloor n/(rk)\rfloor}
f_r\bigl( u_{(i-1)rdk+1},\ldots,u_{irdk}\bigr)\\
% = &\frac{\lfloor n/(rk)\rfloor}{n}\frac1{\lfloor n/(rk)\rfloor}
%f_r( u_{(i-1)rdk+1},\ldots,u_{irdk})\\
&&\quad\to\quad \frac{(r-1)k+1}{rk}v
>\frac{r-1}r v,
\end{eqnarray*}
after using (\ref{eq:overlapornot}).
Taking $r$ as large as we like, we get
${ \lim\inf_{n\to\infty}}
\frac1n\sum_{i=1}^n1_{[\bszero,\bsc)}(\bsx_i)\ge v$.
It follows that
${ \lim\inf_{n\to\infty}}
\frac1n\sum_{i=1}^n1_{[\bsa,\bsb)}(\bsx_i)\ge\vol[\bsa,\bsb)$
for any rectangular subset $[\bsa,\bsb)\subset[0,1]^{dk}$.
Therefore,
${ \lim\sup_{n\to\infty}}
\frac1n\sum_{i=1}^n1_{[\bszero,\bsc)}(\bsx_i)\le v$ too,
for otherwise some rectangle $[\bsa,\bsb)$ would
get too few points.
\end{pf}

%Theorem \ref{thm:ifmAlsFuh} gives sufficient conditions

Now, we prove the main theorems
from Section \ref{sec:consistency}.
\begin{pf*}{ Proof of Theorem \ref{thm:withhomecouple}}
Pick $\varepsilon>0$.
Now let $m\in\Natural$ and
for $i=1,\ldots,n$ define the sequence
$\bsx'_{i,m,0},\ldots,\bsx'_{i,m,m}\in\xset$
%$\bsx'_{i,m,0},\bsx'_{i,m,1},\ldots,\bsx'_{i,m,m}\in\xset$
%by $\bsx'_{i,m,0} =\psi(\bsu_i)$
%and $\bsx'_{i,m,j} = \phi(\bsx'_{i,m,j-1},\bsu_{i+j})$ for $j=1,
as the Rosenblatt--Chentsov transformation
of $\bsu_i,\ldots,\bsu_{i+m}$.

Suppose that $\phi$ is regular and for a bounded rectangle $[\bsa
,\bsb]\subset\Real^s$,
%% xxx should we use $\Omega$? Presumably it has to be Jordan meas
let $f(\bsx) = 1_{[\bsa,\bsb]}(\bsx)$.
Then
%
%e20 ###
%
\begin{equation}\label{three_epsilon}
\int f(\bsx) \pi(\bsx) \,\rd\bsx- \frac{1}{n} \sum_{i=1}^nf(\bsx_i)
%1_{\bsx_i \in R_{\delta,\bsz}}
= \Sigma_1 + \Sigma_2 + \Sigma_3,
\end{equation}
where
\begin{eqnarray*}
\Sigma_1 & = & \int f(\bsx)\pi(\bsx) \,\rd\bsx- \frac{1}{n} \sum
_{i=1}^n f(\bsx'_{i,m,m}),\\
\Sigma_2 & = & \frac{1}{n} \sum_{i=1}^n f(\bsx'_{i,m,m}) - f(\bsx
_{i+m})
\end{eqnarray*}
and
\[
\Sigma_3  = \frac{1}{n} \sum_{i=1}^n f(\bsx_{i+m}) - f(\bsx_i).
\]

For $\Sigma_1$, notice that $\bsx'_{i,m,m} \in[\bsa,\bsb]$ if
and only if $(v_{d(i-1)+1},\ldots, v_{d(i+m)})$ lies in a
$d(m+1)$-dimensional region $\Bcal_1$. The region $\Bcal_1$ has
volume $\int_{[\bsa,\bsb]} \pi(\bsx) \,\rd\bsx$
because $\Pr( \bsx'_{i,m,m}\in[\bsa,\bsb])$ is
$\int_{[\bsa,\bsb]} \pi(\bsx) \,\rd\bsx$ when $(v_{d(i-1)+1},\ldots
,v_{d(i+m)})\sim U[0,1]^{d(m+1)}$.
It has a Riemann integrable indicator function by hypothesis.
Then because $(v_i)_{i\ge1}$ are CUD,
and using Lemma \ref{lem:tech} with $k=m+1$, we get
\[
|\Sigma_1| = \Biggl| \int f(\bsx) \pi(\bsx) \,\rd\bsx- \frac{1}{n}
\sum_{i=1}^n f(\bsx'_{i,m,m})\Biggr|
\longrightarrow0\qquad \mbox{as } n \rightarrow\infty.
\]

Now, consider $\Sigma_2$. The only nonzero terms arise
when $\bsx_{i+m}\ne\bsx'_{i,m,m}$.
This in turn requires that the coupling region $\Ccal$ %leading to the
%home coupling region
is avoided $m$ consecutive times, by
$\bsu_{i+1},\ldots,\bsu_{i+m}$.
Then
$(v_{di+1},\ldots,v_{d(i+m)})$ belongs to a region
of volume at most $(1-\vol(\Ccal))^m$.
Choose $m$ large enough that $(1-\vol(\Ccal))^m<\varepsilon$.
Then
\[
\lim\sup_{n\to\infty}
\Biggl|\frac{1}{n} \sum_{i=1}^n f(\bsx'_{i,m,m})- f(\bsx
_{i+m})\Biggr|
< \varepsilon.
\]

For the third term, $|\Sigma_3|$ is at most $m/n$, which goes to $0$
as $n
\rightarrow\infty$.
Thus, we have
\[
\Biggl|\lim_{n\rightarrow\infty} \frac{1}{n} \sum_{i=1}^n 1_{\bsx
_i\in[\bsa,\bsb]}
- \int_{[\bsa,\bsb]} \pi(\bsx) \,\rd\bsx
\Biggr| < \varepsilon.
\]
As $\varepsilon> 0$ was chosen arbitrarily, the result follows for
this case.

The result holds trivially for the function $1_{\xset}$, hence we are done.
% If instead the MCMC is regular for bounded continuous functions, then
%we still
% make the decomposition \eqref{three_epsilon}, except that now $f$
% is a bounded continuous function on $\Real^s$.
% The three terms converge as before and so $\bsx_1,\ldots,\bsx_n$
% consistently samples $\pi$.% by Proposition \ref{prop:forcts}.
\end{pf*}
\begin{pf*}{Proof of Theorem \ref{thm:ifmAlsFuh}}
We use the notation from the proof of
Theorem~\ref{thm:withhomecouple}. As in the proof of
Theorem \ref{thm:withhomecouple}, we write $\int f(\bsx)\pi(\bsx
)\,\rd
\bsx-\frac1n\sum_{i=1}^n f(\bsx_i)$ as the sum of three terms. The
first and third terms vanish by the same arguments we used in
Theorem \ref{thm:withhomecouple}.

For the second term, we have
\[
|\Sigma_2(n)| \le\frac{1}{n} \sum_{i=1}^n |f(\bsx'_{i,m,m}) -
f(\bsx_{i+m})|.
\]
Let $\varepsilon> 0$ be arbitrary. We show that
$\limsup_{n\rightarrow\infty} |\Sigma_2(n)| \le\varepsilon$. As
$\varepsilon> 0$ is arbitrary, this then implies that $\limsup_{n
\rightarrow\infty} |\Sigma_2(n)| = 0$.

Assume that the Gibbs sampler is regular for rectangles and
for a bounded positive volume rectangle $[\bsa,\bsb] \subset\Real^s$
let $f(\bsx) = 1_{[\bsa,\bsb]}(\bsx)$.
For $0\le\delta< \break\min_{1\le j\le d}(b_j-a_j)$, let
$\bsdelta= (\delta,\ldots,\delta)\in\Real^s$ and put
$f_\delta(\bsx) = 1_{[\bsa-\bsdelta,\bsb+\bsdelta]}$
and $f_{-\delta}(\bsx) = 1_{[\bsa+\bsdelta,\bsb-\bsdelta]}$.

Because $f_\delta(\bsx)\ge f(\bsx)\ge f_{-\delta}(\bsx)$,
the triple $(f_{-\delta}(\bsx'_{i,m,m}),f(\bsx'_{i,m,m}),\break
f_{\delta}(\bsx'_{i,m,m}))$
must be in the set $S=\{ (0,0,0),(0,0,1),(0,1,1),(1,1,1)\}$.
Likewise $f(\bsx_{i+m})\in\{0,1\}$.
By inspecting all $8$ cases in $S\times\{0,1\}$, we
find that
$|\Sigma_2|\le\sigma_1+\sigma_2+\sigma_3$, for
\begin{eqnarray*}
\sigma_1 & = & \frac{1}{n} \sum_{i=1}^n f_\delta(\bsx'_{i,m,m}) -
f_{-\delta}(\bsx'_{i,m,m}),\\
\sigma_2 & = & \frac{1}{n} \sum_{i=1}^n \bigl(f_{-\delta}(\bsx
'_{i,m,m})-f(\bsx_{i+m})\bigr)_+
\end{eqnarray*}
and
\[
\sigma_3 = \frac{1}{n} \sum_{i=1}^n \bigl(f(\bsx_{i+m})-f_{\delta
}(\bsx'_{i,m,m})\bigr)_+,
\]
where $z_+=\max(z,0)$.

Choose $\delta> 0$ such that
\[
\int_{\xset\cap
([\bsa-\bsdelta,\bsb+\bsdelta] \setminus[\bsa+\bsdelta,\bsb
-\bsdelta])} \pi(\bsx) \,\rd
\bsx< \frac{\varepsilon}{3}.% \mbox{ and } \int_{\xset\cap
%([\bsa,\bsb] \setminus[\bsa+\bsdelta,\bsb-\bsdelta])} \pi(\bsx) \,\rd
\]
As the Gibbs sampler is regular for rectangles, $(v_i)_{i\ge1}$ is a
CUD sequence, and $\bsx'_{i,m,m}$ is constructed using the
Rosenblatt--Chentsov transformation we have
\begin{eqnarray*}
& & \lambda(\{\bsu\in[0,1]^{dm+d}\dvtx \bsx'_{i,m,m} \in
[\bsa-\bsdelta,\bsb+\bsdelta]\setminus[\bsa,\bsb] \}
) \\
&&\qquad = \int_{\xset\cap
([\bsa-\bsdelta,\bsb+\bsdelta] \setminus[\bsa+\bsdelta,\bsb
+\bsdelta])} \pi(\bsx) \,\rd
\bsx< \frac{\varepsilon}{3},
\end{eqnarray*}
and so $\limsup_{n\rightarrow\infty} |\sigma_1(n)| \le{\varepsilon}/{3}$.
%Analogously we obtain $\limsup_{n\rightarrow\infty} |\sigma_4(n)| \le

The points $\bsx'_{i,m,m}$ and $\bsx_{i+m}$ have different starting
points $\bsx'_{i,m,0}$ and $\bsx_{i}$, but are updated $m$ times
using the same $\bsu_{i+1},\ldots, \bsu_{i+m}$, that is, $\bsx
'_{i,m,m} = \phi_m(\bsx'_{i,m,0}$, $\bsu_{i+1},\ldots,\bsu_m)$ and
$\bsx_{i+m} = \phi_m(\bsx_{i}, \bsu_{i+1},\ldots, \bsu_{i+m})$.
Therefore, Theorem \ref{thm:ifmAlsFuh} implies that there is a
constant $C > 0$ such that for all sufficiently large $m \ge m_i^\ast$
the region % xxx Josef: I made m^* into m_i^*
\begin{eqnarray*}
\Bcal_{m,i} & = & \{(\bsv_1,\ldots, \bsv_m) \in[0,1]^{dm}\dvtx
d(\phi_m(\bsx'_{i,m,0},(\bsv_1,\ldots,\bsv_m)),\\
&&\hspace*{142.2pt}\phi_m(\bsx
_{i},(\bsv_1,\ldots,\bsv_m))) > \gamma^m\},
\end{eqnarray*}
has volume at most $C \alpha_\gamma^m$. Let $\Bcal_m =
\bigcup_{i=1}^n \Bcal_{m,i}$. Let $\beta= \infty$ if $[\bsa,\bsb]
\cap\xset= \varnothing$ or $\xset\setminus[\bsa-\bsdelta,\bsb
+\bsdelta] = \varnothing$ and $\beta= \inf\{d(\bsy,\bsy')\dvtx \bsy\in
[\bsa,\bsb] \cap\xset, \bsy' \in\xset\setminus[\bsa-\bsdelta
,\bsb+\bsdelta]\}$ otherwise.

%Let $m_1 = m_1(n) > m^\ast$ be such % xxx I changed this
Let $m_1 = m_1(n)$ be such
that $C n \alpha_\gamma^{m_1} < \varepsilon/3$ and $\gamma^{m_1}
< \beta$. Now take $m_0 \ge\max\{m_1,m_1^\ast,\ldots,m_n^\ast\}$.
For large enough $n$, we
can take $m_0 = m_0(n) = \lceil\frac{\log n + \log
(2C/\varepsilon)}{\log1/\alpha_\gamma}\rceil+1$. Then
% xxx Josef: I moved /\varepsilon into the log( ) function
$\Bcal_{m_0}$ has volume at most $\varepsilon/3$.\vspace*{1pt}

Thus,\vspace*{1pt} $f_{-\bsdelta}(\bsx'_{i,m_0,m_0}) >
f(\bsx_{i+m_0})$ implies that
$d(\bsx'_{i,m_0,m_0},\bsx_{i+m_0}) \ge\beta$, which in turn implies
that $(\bsu_{i+1},\ldots, \bsu_{i+m_0}) \in\Bcal_{m_0,i}$, and
%therefore also % therefore also ... gives a bad line break here
so
$(\bsu_{i+1},\ldots, \bsu_{i+m_0}) \in\Bcal_{m_0}$. Therefore, we have
\[
{\limsup_{n\rightarrow\infty} }|\sigma_2(n)| \le\limsup
_{n\rightarrow\infty} \frac{1}{n} \sum_{i=1}^n 1_{(\bsu
_{i+1},\ldots, \bsu_{i+m_0}) \in\Bcal_{m_0}} = \limsup_{m_0
\rightarrow\infty} \lambda(\Bcal_{m_0}) \le\frac{\varepsilon}{3}.
\]
A similar argument shows that $\limsup_{n\rightarrow\infty}|\sigma
_3(n)|\le\varepsilon/3$.

Combining the three bounds yields
\begin{eqnarray*}
{\limsup_{n\rightarrow\infty} |\Sigma_2(n)| } & \le& \limsup
_{n\rightarrow\infty} \sigma_1(n) +\limsup_{n\rightarrow\infty}
\sigma_2(n) +\limsup_{n\rightarrow\infty} \sigma_3(n) \\ & \le&
\frac{\varepsilon}{3} + \frac{\varepsilon}{3} + \frac{\varepsilon
}{3} = \varepsilon,
\end{eqnarray*}
establishing consistency when the Gibbs sampler is regular.

Since the result holds trivially for the function $1_{\xset}$, the
result follows.
\end{pf*}

The coupling region in Theorem \ref{thm:withhomecouple} was replaced
by a mean contraction assumption $\int_{[0,1]^d} \log(\ell(\bsu
))\,\rd\bsu< 0$ in Theorem \ref{thm:withiterfunction}.\vspace*{1pt} This way we
obtain (possibly different) coupling type regions $\Bcal_{m,i}$ for
each $i = 1,\ldots, n$. We remedy this situation by letting $m$ depend
on $n$, which in turn requires us to use a stronger assumption on the
CUD sequence $(v_i)_{i\ge1}$, namely, that $\lim_{n\rightarrow\infty
} D^{\ast d_n}_n = 0$.
\end{appendix}

\section*{Acknowledgments}

We thank Seth Tribble, Erich Novak, Ilya M. Sobol' and two
anonymous reviewers for helpful comments.

%suskaldyti doi

%
\printaddresses


\begin{thebibliography}{99}

%b1 ###
\bibitem{albechib1993}
\textsc{Albert}, J. and \textsc{Chib}, S. (1993).
{Bayesian} analysis of binary and polychotomous response data.
\textit{J. Amer. Statist. Assoc.} \textbf{88} 669--679.
\MR{1224394}


%b2 ###
\bibitem{alsmfuh2001}
\textsc{Alsmeyer}, G. and \textsc{Fuh}, C.-D. (2001).
Limit theorems for iterated function mappings.
\textit{Stochastic Process. Appl.} \textbf{96} 123--142.
\MR{1856683}

%b3 ###
\bibitem{andrmoul2006}
\textsc{Andrieu}, C. and \textsc{Moulines}, E. (2006).
On the ergodicity properties of some {MCMC} algorithms.
\textit{Ann. Appl. Probab.} \textbf{16} 1462--1505.
\MR{2260070}

%b4 ###
\bibitem{ash1972}
\textsc{Ash}, R. B. (1972).
\textit{Real Analysis and Probability}.
Academic Press, New York.
\MR{0435320}

%b5 ###
\bibitem{bill1999}
\textsc{Billingsley}, P. (1999).
\textit{Convergence of Probability Measures}.
Wiley, New York.
\MR{1700749}

%b6 ###
\bibitem{chau2004}
\textsc{Chaudary}, S. (2004).
Acceleration of {Monte Carlo} methods using low discrepancy
sequences.
Ph.D. thesis, UCLA.

%b7 ###
\bibitem{chen1967}
\textsc{Chentsov}, N. N. (1967).
Pseudorandom numbers for modelling {Markov} chains.
\textit{Comput. Math. Math. Phys.} \textbf{7} 218--2332.

%b8 ###
\bibitem{crailemi2007}
\textsc{Craiu}, R. V. and \textsc{Lemieux}, C. (2007).
Acceleration of the multiple-try {Metropolis} algorithm using
antithetic and stratified sampling.
\textit{Stat. Comput.} \textbf{17} 109--120.
\MR{2380640}

%b9 ###
\bibitem{devr1986}
\textsc{Devroye}, L. (1986).
\textit{Non-Uniform Random Variate Generation}.
Springer, New York.
\MR{0836973}

%b10 ###
\bibitem{diacfree1999}
\textsc{Diaconis}, P. and \textsc{Freedman}, D. (1999).
Iterated random functions.
\textit{SIAM Rev.} \textbf{41} 45--76.
\MR{1669737}

%b11 ###
\bibitem{dicknote2007}
\textsc{Dick}, J. (2007).
A note on the existence of sequences with small star discrepancy.
\textit{J. Complexity} \textbf{23} 649--652.
\MR{2372019}

%b12 ###
\bibitem{dickmcqmc2009}
\textsc{Dick}, J. (2009).
On {quasi-Monte Carlo} rules achieving higher order convergence.
In \textit{{Monte Carlo} and
Quasi-{Monte Carlo} Methods 2008} (P. L'Ecuyer and A. B. Owen, eds.)
73--96.
Springer, Heidelberg.

%b13 ###
\bibitem{doerfrie2009}
\textsc{Doerr}, B. and \textsc{Friedrich}, T. (2009).
Deterministic random walks on the two-dimensional grid.
\textit{Combin. Probab. Comput.} \textbf{18} 123--144.
\MR{2497377}

%b14 ###
\bibitem{finn1947}
\textsc{Finney}, D. J. (1947).
The estimation from individual records of the relationship between
dose and quantal response.
\textit{Biometrika} \textbf{34} 320--334.

%b15 ###
\bibitem{gaveomui1987}
\textsc{Gaver}, D. and \textsc{O'Murcheartaigh}, I. (1987).
Robust empirical {Bayes} analysis of event rates.
\textit{Technometrics} \textbf{29} 1--15.
\MR{0876882}

%b16 ###
\bibitem{gelfsmit1990}
\textsc{Gelfand}, A. E. and \textsc{Smith}, A. F. M. (1990).
Sampling-based approaches to calculating marginal densities.
\textit{J. Amer. Statist. Assoc.} \textbf{85} 398--409.
\MR{1141740}

%b17 ###
\bibitem{gelmshir2010}
\textsc{Gelman}, A. and \textsc{Shirley}, K. (2010).
Inference from simulations and monitoring convergence.
In \textit{Handbook of {Markov Chain Monte Carlo}: Methods and Applications}.
(S.~Brooks, A. Gelman, G. Jones and X.-L. Meng, eds.) 131--143.
Chapman and Hall/CRC Press, Boca Raton, FL.

%b18 ###
\bibitem{ghorlima2006}
\textsc{Ghorpade}, S. R. and \textsc{Limaye}, B. V. (2006).
\textit{A Course in Calculus and Real Analysis}.
Springer, New York.
\MR{2229667}

%b19 ###
\bibitem{gnewsrivwinz2008}
\textsc{Gnewuch}, M., \textsc{Srivastav}, A. and \textsc{Winzen}, C. (2008).
Finding optimal volume subintervals with $k$ points and
computing the
star discrepancy are {NP}-hard.
\textit{J. Complexity} \textbf{24} 154--172.

%b20 ###
\bibitem{gord1941}
\textsc{Gordon}, R. D. (1941).
Value of {Mill's} ratio of area to bounding ordinate and of the
normal probability integral for large values of the argument.
\textit{Ann. Math. Statist.} \textbf{18} 364--366.
\MR{0005558}

%b21 ###
\bibitem{haarsakstamm2001}
\textsc{Haario}, H., \textsc{Saksman}, E. and \textsc{Tamminen}, J. (2001).
An adaptive {Metropolis} algorithm.
\textit{Bernoulli} \textbf{7} 223--242.
\MR{1828504}

%b22 ###
\bibitem{knut199723}
\textsc{Knuth}, D. E. (1998).
\textit{The Art of Computer Programming}, 3rd ed.
\textit{Seminumerical Algorithms}~\textbf{2}.
Addison-Wesley, Reading, MA.
\MR{0378456}

%b23 ###
\bibitem{lebe1902}
\textsc{Lebesgue}, H. L. (1902).
Int\'egrale, longueur, aire.
Ph.D. thesis, Univ. de Paris.

%b24 ###
\bibitem{leculecotuff2008}
\textsc{L'Ecuyer}, P., \textsc{Lecot}, C. and \textsc{Tuffin}, B. (2008).
A randomized {quasi-Monte Carlo} simulation method for {Markov}
chains.
\textit{Oper. Res.} \textbf{56} 958--975.

%b25 ###
\bibitem{leculemi1999a}
\textsc{L'Ecuyer}, P. and \textsc{Lemieux}, C. (1999).
Quasi-{Monte Carlo} via linear shift-register sequences.
In \textit{Proceedings of the 1999 Winter Simulation Conference}
(P. A. Farrington, H. B. Nembhard, D. T. Sturrock and
G. W. Evans, eds.) 632--639. IEEE Press, Piscataway, NJ.

%b26 ###
\bibitem{lemisido2006}
\textsc{Lemieux}, C. and \textsc{Sidorsky}, P. (2006).
Exact sampling with highly uniform point sets.
\textit{Math. Comput. Modelling} \textbf{43} 339--349.
\MR{2214643}

%b27 ###
\bibitem{liao1998}
\textsc{Liao}, L. G. (1998).
Variance reduction in {Gibbs} sampler using quasi random numbers.
\textit{J. Comput. Graph. Statist.} \textbf{7} 253--266.

%b28 ###
\bibitem{liu2001}
\textsc{Liu}, J. S. (2001).
\textit{{Monte Carlo} Strategies in Scientific Computing}.
Springer, New York.
\MR{1842342}

%b29 ###
\bibitem{matsnish1998}
\textsc{Matsumoto}, M. and \textsc{Nishimura}, T. (1998).
{Mersenne} twister: A 623-dimensionally equidistributed uniform
pseudorandom number generator.
\textit{ACM Transactions on Modeling and Computer Simulation} \textbf{8} 3--30.

\bibitem{marshoff1993}
\textsc{Marsden, J. E.} and \textsc{Hoffman}, M. J.
(1993). \textit{Elementary Classical Analysis}, 2nd ed.
Macmillan, New York.

%b30 ###
\bibitem{morocafl1993}
\textsc{Morokoff}, W.  and \textsc{Caflisch}, R. E. (1993).
A quasi-{Monte Carlo} approach to particle simulation of the
heat equ
ation.
\textit{SIAM J. Numer. Anal.} \textbf{30} 1558--1573.
\MR{1249033}

%b31 ###
\bibitem{neal2003}
\textsc{Neal}, R. M. (2003).
Slice sampling.
\textit{Ann. Statist.} \textbf{31} 705--767.
\MR{1994729}

%b32 ###
\bibitem{nied1986}
\textsc{Niederreiter}, H. (1986).
Multidimensional integration using pseudo-random numbers.
\textit{Math. Programming Stud.} \textbf{27} 17--38.
\MR{0836749}

%b33 ###
\bibitem{nied92}
\textsc{Niederreiter}, H. (1992).
\textit{Random Number Generation and Quasi-{Monte Carlo} Methods}.
SIAM, Philadelphia, PA.
\MR{1172997}

%b34 ###
\bibitem{rtms}
\textsc{Owen}, A. B. (1995).
Randomly permuted $(t,m,s)$-nets and $(t,s)$-sequences.
In \textit{Monte Carlo and Quasi-Monte Carlo Methods in Scientific Computing}
(H. Niederreiter and P. Jau-Shyong Shiue, eds.) 299--317.
Springer, New York.
\MR{1445791}

%b35 ###
\bibitem{variation}
\textsc{Owen}, A. B. (2005).
Multidimensional variation for quasi-{Monte Carlo}.
In \textit{Contemporary Multivariate Analysis and Design
of Experiments: In Celebration of
Prof. Kai-Tai Fang's 65th Birthday}
(J. Fan and G. Li, eds.).
World Sci. Publ., Hackensack, NJ.
\MR{2271076}

%b36 ###
\bibitem{qmcmetro}
\textsc{Owen}, A. B. and \textsc{Tribble}, S. D. (2005).
A quasi-{Monte Carlo} {Metropolis} algorithm.
\textit{Proc. Natl. Acad. Sci. USA} \textbf{102} 8844--8849.
\MR{2168266}

%b37 ###
\bibitem{propwils1996}
\textsc{Propp}, J. G. and \textsc{Wilson}, D. B. (1996).
Exact sampling with coupled {Markov} chains.
\textit{Random Structures and Algorithms} \textbf{9} 223--252.
\MR{1611693}

%b38 ###
\bibitem{robecase2004}
\textsc{Robert}, C. P. and \textsc{Casella}, G. (2004).
\textit{Monte Carlo Statistical Methods}, 2nd ed.
Springer, New York.
\MR{2080278}

%b39 ###
\bibitem{roberoseschw1998}
\textsc{Roberts}, G. O., \textsc{Rosenthal}, J. S. and \textsc{Schwartz}, P. O. (1998).
Convergence properties of perturbed {Markov} chains.
\textit{J. Appl. Probab.} \textbf{35} 1--11.
\MR{1622440}

%b40 ###
\bibitem{rose1952}
\textsc{Rosenblatt}, M. (1952).
Remarks on a multivariate transformation.
\textit{Ann. Math. Statist.} \textbf{23} 470--472.
\MR{0049525}

%b41 ###
\bibitem{sobo1974}
\textsc{Sobol'}, I. M. (1974).
Pseudo-random numbers for constructing discrete {Markov}
chains by
the {Monte Carlo} method.
\textit{USSR Comput. Math. Math. Phys.} \textbf{14} 36--45.
\MR{0339444}

%b42 ###
\bibitem{trib2007}
\textsc{Tribble}, S. D. (2007).
Markov chain Monte Carlo algorithms using completely uniformly
distributed driving sequences.
Ph.D. thesis, Stanford Univ.
\MR{2710331}

%b43 ###
\bibitem{qmcmetro2}
\textsc{Tribble}, S. D. and \textsc{Owen}, A. B. (2008).
Construction of weakly {CUD} sequences for {MCMC} sampling.
\textit{Electron. J. Stat.} \textbf{2} 634--660.
\MR{2426105}

%b44 ###
\bibitem{weyl1916}
\textsc{Weyl}, H. (1916).
\"{U}ber die gleichverteilung von zahlen mod. eins.
\textit{Math. Ann.} \textbf{77} 313--352.
\MR{1511862}

\end{thebibliography}
\end{document}